\documentclass[12pt]{article}
\usepackage{amsmath,amssymb,bm}
\usepackage[dvips]{graphics}
\usepackage{color}

\newtheorem{thm}{Theorem}[subsection]
\newtheorem{lem}[thm]{Lemma}
\newtheorem{prop}[thm]{Proposition}
\newtheorem{cor}[thm]{Corollary}
\newtheorem{define}[thm]{Definition}
\newtheorem{ex}[thm]{Example}
\newtheorem{rem}[thm]{Remark}
\newtheorem{claim}{Claim}
\newtheorem{THM}{Theorem}

\newcommand{\qed}{\hfill$\Box$\par}
\newcommand{\bqed}{\hfill$\blacksquare$\par}

\begin{document}
\baselineskip=20pt
\begin{center}
\textbf{\Large Demazure Character Formulas for Generalized Kac--Moody Algebras}
\end{center}
\begin{center}
{\large Motohiro Ishii}\\ 
Graduate School of Pure and Applied Sciences, University of Tsukuba, \\
Tsukuba, Ibaraki 305-8571, Japan \\
(e-mail: \verb/ishii731@math.tsukuba.ac.jp/)
\end{center}
\begin{quote}
\begin{center}
\textbf{Abstract}
\end{center}
For a dominant integral weight $\lambda $,
we introduce a family of $U_q ^+ (\mathfrak{g})$-submodules $V_w (\lambda )$
of the irreducible highest weight $U_q (\mathfrak{g})$-module $V(\lambda )$ 
of highest weight $\lambda $ for a
generalized Kac--Moody algebra $\mathfrak{g}$. 
We prove that the module $V_w (\lambda )$ is spanned by its 
global basis, and then give a character formula for $V_w (\lambda )$, 
which generalizes the Demazure character formula
for ordinary Kac--Moody algebras.
\end{quote}
\section{Introduction}
For a (symmetrizable) Kac--Moody algebra $\mathfrak{g}$, 
the \textit{Demazure\ character\ formula} describes the (formal)
character of a \textit{Demazure\ module}, which is a $U^+ (\mathfrak{g})$-submodule 
generated by an extremal weight vector of an integrable 
highest weight $U(\mathfrak{g})$-module; this formula was proved by 
Kumar ([\textbf{Kum}]) and Mathieu ([\textbf{M}])
independently by using geometric methods.
Also, from an algebraic viewpoint, it is possible to formulate the notion 
of Demazure modules for an integrable highest weight $U_q (\mathfrak{g})$-module. 
In fact, Littelmann ([\textbf{Li3}]) gave a (conjectural) algebraic description of the Demazure 
character formula in terms of Kashiwara's crystal bases for $\mathfrak{g}$ of finite type. 
Soon afterward, this conjecture of Littelmann was proved by Kashiwara ([\textbf{Kas2}])
generally for a Kac--Moody algebra. More precisely, Kashiwara showed that a 
Demazure module is spanned by its global basis, and that there exists 
a basis at the crystal limit $``q=0"$, which is called a \textit{Demazure\ crystal}.

In [\textbf{JKK,\ JKKS}], Kashiwara's crystal basis theory was extended to 
the case of generalized Kac--Moody algebras. Also, Littelmann's path model for representations of 
Kac--Moody algebras was extended to the case of generalized Kac--Moody algebras
by Joseph and Lamprou ([\textbf{JL}]). Therefore, it is natural to expect the existence of 
Demazure crystals for an irreducible highest weight module over a generalized Kac--Moody algebra.
The purpose of this paper is to obtain a generalization of Demazure modules, their crystal bases, 
and a character formula for them for generalized Kac--Moody algebras.

It is known that the representation theory of generalized Kac--Moody algebras
are very similar to that of ordinary Kac--Moody algebras, and many 
results for Kac--Moody algebras are extended to the case of generalized 
Kac--Moody algebras. However, there are some obstructions, 
which come from the existence of \textit{imaginary\ simple\ roots},  for the study of the structure 
of an irreducible highest weight module over a generalized Kac--Moody algebra.
Recall that a Demazure crystal decomposes into a disjoint union of 
$\mathfrak{sl}_2$-strings of finite length in the case of ordinary Kac--Moody algebras. 
In contrast, in the case of generalized Kac--Moody algebras, 
$\mathfrak{sl}_2$-strings and ``Heisenberg algebra"-strings 
corresponding to imaginary simple roots for an irreducible highest weight 
module are no longer of finite length.
Hence, we cannot apply a method similar to the case of ordinary Kac--Moody algebras to our setting. 
To overcome this difficulty, we introduce a certain Kac--Moody algebra $\Tilde{\mathfrak{g}}$ 
defined from a given Borcherds--Cartan datum of a generalized Kac--Moody algebra $\mathfrak{g}$, 
and then we relate the representation theory of $\mathfrak{g}$ to the one of $\Tilde{\mathfrak{g}}$. 
This enables us to study the structure of an irreducible highest weight module over a 
generalized Kac--Moody algebra by comparing it with the corresponding one 
over an ordinary Kac--Moody algebra. In this way, under a certain condition, 
we can define Demazure modules for a generalized Kac--Moody algebra, 
and show that they are spanned by their global bases (Theorem \ref{thm3}). 
As a result, we obtain a character formula for these Demazure modules 
for a generalized Kac--Moody algebra by introducing a ``modified" Demazure operator 
for each imaginary simple root (Theorem \ref{thm4}).

Let us state our results more precisely. Let $\mathfrak{g}$
be a generalized Kac--Moody algebra and $\Tilde{\mathfrak{g}}$ 
the associated Kac--Moody algebra  (see \S 2.2).
For a dominant integral weight $\lambda \in P^+$ for $\mathfrak{g}$ 
and the corresponding dominant integral weight 
$\Tilde{\lambda } \in \widetilde{P}^+$ for $\Tilde{\mathfrak{g}}$ (see \S 2.5),
we denote by $\mathbb{B}(\lambda )$ and 
$\widetilde{\mathbb{B}}(\Tilde{\lambda })$ the sets of 
\textit{generalized Lakshmibai--Seshadri paths} and 
ordinary \textit{Lakshmibai--Seshadri paths}, respectively.
Then, there exists an embedding 
$\mathbb{B}(\lambda ) \hookrightarrow 
\widetilde{\mathbb{B}}(\Tilde{\lambda })$ 
of path crystals (Proposition \ref{2.5.1});
note that this is not a morphism of crystals.
By using this embedding, we have the following 
decomposition rules for $U_q (\mathfrak{g})$-modules 
in the category $\mathcal{O}_{\mathrm{int}}$ (see [\textbf{JKK}, \S 2]).
Denote by $V(\lambda )$ the irreducible highest weight 
$U_q (\mathfrak{g})$-module of highest weight $\lambda \in P^+$.
\begin{THM}\label{thm1}
Let $\lambda ,\mu \in P^+$ be dominant integral weights for $\mathfrak{g}$.
Then, we have an isomorphism of $U_q (\mathfrak{g})$-modules:
$$V (\lambda ) \otimes V (\mu )\ \cong \ 
\bigoplus _
{\begin{subarray}{c}
\pi \in \mathbb{B}(\mu ) \\ 
\Tilde{\pi }\ \!\!:\ \!\! \Tilde{\lambda }
\text{-}\mathrm{dominant}
\end{subarray}}
V \bigl( \lambda +\pi (1) \bigl).$$
Here, $\Tilde{\pi } \in \widetilde{\mathbb{B}}(\Tilde{\mu })$ denotes the 
image of $\pi \in \mathbb{B}(\mu )$ under the embedding 
$\mathbb{B}(\mu ) \hookrightarrow \widetilde{\mathbb{B}}(\Tilde{\mu })$,
and it is said to be $\Tilde{\lambda }$-dominant
if $\Tilde{\pi }(t)+\Tilde{\lambda }$ belongs to the dominant Weyl chamber 
of $\Tilde{\mathfrak{g}}$ for all $t\in [0,1]$.
\end{THM}
\begin{THM}\label{thm2}
Let $\lambda \in P^+$. For a subset $S$ of the index set $I$ of 
simple roots of $\mathfrak{g}$, we denote by $\mathfrak{g}_S$ the corresponding Levi subalgebra 
of $\mathfrak{g}$, by $\Tilde{\mathfrak{g}}_{\widetilde{S}}$ the corresponding 
one for $\Tilde{\mathfrak{g}}$ (see \S 3.1), and by $V_S (\mu )$ the irreducible highest weight 
$U_q (\mathfrak{g}_S)$-module of highest weight $\mu $. Then,
$$V(\lambda )\ \cong \ \bigoplus _
{\begin{subarray}{c}\pi \in \mathbb{B}(\lambda) \\ 
\Tilde{\pi }\ \!\! :\ \!\! \Tilde{\mathfrak{g}}_{\widetilde{S}}
\text{-}\mathrm{dominant}
\end{subarray}} V_S \bigl( \pi (1) \bigl) \ \ 
\mathrm{as}\ U_q (\mathfrak{g}_S ) \mathrm{\text{-}modules}.$$
Here, $\Tilde{\pi } \in \widetilde{\mathbb{B}}(\Tilde{\lambda })$ 
denotes the image of $\pi \in \mathbb{B}(\lambda )$ under the 
embedding $\mathbb{B}(\lambda ) \hookrightarrow \widetilde{\mathbb{B}}(\Tilde{\lambda })$,
and it is said to be $\Tilde{\mathfrak{g}}_{\widetilde{S}}$-dominant
if $\Tilde{\pi }(t)$ belongs to the dominant Weyl chamber of $\Tilde{\mathfrak{g}}_{\widetilde{S}}$
for all $t\in [0,1]$.
\end{THM}
Also, we prove an analog of the  
\textit{Parthasarathy--Ranga Rao--Varadarajan conjecture} 
for generalized Kac--Moody algebras (Theorem \ref{3.2.4}), using 
the embedding of path models and the tensor product decomposition rule (Theorem \ref{thm1}).

Now, let us describe our main results. For a dominant integral weight $\lambda \in P^+$, 
we take and fix an element $w$ in the monoid $\mathcal{W}$ (see Definition \ref{2.3.1}) 
satisfying the condition: 
\begin{quote}
there exists an expression 
$w=v_l r_{j_l} \cdots v_1 r_{j_1} v_0$, with $j_1 ,\ldots ,j_l \in I^{im}$ 
and $\ v_0 ,v_1 ,\ldots ,v_l \in \mathcal{W}_{re}$ (see \S 2.1 and \S 2.3), such that
\begin{align*}
\alpha _{j_s} ^{\vee}\bigl( v_{s-1} r_{j_{s-1}} \cdots 
v_1 r_{j_1} v_0 (\lambda )\bigl) =1\ \mathrm{for\ all}\ s=1,2,\ldots ,l,
\end{align*}
where $\alpha _j ^{\vee},\ j\in I,$ are the simple coroots of $\mathfrak{g}$. 
\end{quote}
For such an element $w\in \mathcal{W}$, we set
$V_w (\lambda ):=U_q ^+ (\mathfrak{g}) \ \! V(\lambda )_{w\lambda }$, 
where $V(\lambda )_{\mu } \subset V(\lambda )$ 
is the weight space of weight $\mu \in P$.
\begin{THM}\label{thm3}
For $\lambda \in P^+$ and $w\in \mathcal{W}$ satisfying the condition above,
there exists a subset $B_w (\lambda )\subset B(\lambda )$
of the crystal basis $B(\lambda )$ of $V(\lambda )$ such that 
$$V_w (\lambda )=\bigoplus _{b\in B_w (\lambda )}
\mathbb{C}(q)G_{\lambda }(b),$$
where $\{ G_{\lambda }(b) \} _{b\in B(\lambda )}$
denotes the global basis of $V(\lambda )$.
\end{THM}
For the element $w \in \mathcal{W}$ above, 
we can take a specific expression, 
a \textit{minimal dominant reduced expression} (see \S 4.2), 
$w =w_k r_{i_k} ^{a_k} \cdots w_1 r_{i_1 }^{a_1 } w_0$, 
with $w_0 ,w_1 ,\ldots ,w_k \in \mathcal{W}_{re},\ 
a_1 ,\ldots ,a_k \in \mathbb{Z}_{>0}$, and 
$i_1 ,\ldots ,i_k \in I^{im}$ (all distinct).
For each $i\in I^{re}$, we define an operator $\mathcal{D}_i$ 
on the group ring $\mathbb{Z}[P]:=\bigoplus _{\mu \in P} \mathbb{Z} e^{\mu }$
of the weight lattice $P$ of $\mathfrak{g}$ by: 
$$\mathcal{D}_i (e^{\mu }):=
\frac{e^{\mu }-e^{\mu -(1+\alpha _i ^{\vee}(\mu ))\alpha _i}}
{1-e^{-\alpha _i}}.$$
For a reduced expression 
$v=r_{j_l} \cdots r_{j_2} r_{j_1} \in \mathcal{W}_{re}$, with 
$\ j_1 ,j_2 ,\ldots ,j_l \in I^{re}$,
we set $\mathcal{D}_v := \mathcal{D}_{j_l} \cdots \mathcal{D}_{j_2}\mathcal{D}_{j_1}$. 
Also, for each $i\in I^{im}$ and $a\in \mathbb{Z}_{\ge 1}$, 
we define an operator $\mathcal{D}_i ^{(a)}$ by
\begin{align*}
\mathcal{D}_i ^{(a)}(e^{\mu }):=
\begin{cases}
e^{\mu } & \mathrm{if}\ \alpha _i ^{\vee}(\mu )=0, \\
\sum _{m=0} ^a e^{\mu -m\alpha _i} & \mathrm{otherwise}.
\end{cases}
\end{align*}
\begin{THM}\label{thm4}
Let $\lambda \in P^+$ and $w = w_k r_{i_k} ^{a_k} \cdots w_1 r_{i_1 }^{a_1 } w_0 \in \mathcal{W}$ 
be as above. Then, the (formal) character of the Demazure module 
$V_w (\lambda )$ is given as follows:
$$\mathrm{ch}\ \! V_w (\lambda ) 
=\mathcal{D}_{w_k} \mathcal{D}_{i_k} ^{(a_k)} \cdots 
\mathcal{D}_{w_1} \mathcal{D}_{i_1} ^{(a_1)} \mathcal{D}_{w_0} (e^{\lambda }).$$
\end{THM}

This paper is organized as follows. In Section 2, we recall some elementary facts 
about generalized Kac--Moody algebras and Joseph--Lamprou's path model.
Also, we review the construction of an embedding of Joseph--Lamprou's path model into 
Littelmann's path model. In Section 3, we give proofs of Theorems \ref{thm1} and \ref{thm2}, 
and show an analog of the Parthasarathy--Ranga-Rao--Varadarajan conjecture 
for generalized Kac--Moody algebras. 
In Section 4, we introduce Demazure modules for 
generalized Kac--Moody algebras, and study the structure of them. 
Finally in Section 5, we prove Theorems \ref{thm3} and \ref{thm4}.
\begin{flushleft}
\textbf{Acknowledgments}
\end{flushleft}
The author is grateful to Professor Satoshi Naito 
for reading the manuscript very carefully and for valuable comments. 
He is also grateful to Professor Daisuke Sagaki for valuable comments. 
The author's research was supported by the Japan Society for the 
Promotion of Science Research Fellowships for Young Scientists.

\section{Preliminaries}
\subsection{Generalized Kac--Moody algebras}
In this subsection, we recall some fundamental facts about 
generalized Kac--Moody algebras. For more details, 
we refer the reader to [\textbf{Bo}, \textbf{JKK}, 
\textbf{JKKS}, \textbf{JL}, \textbf{Kac}, \textbf{Kan}].

Let $I$ be a (finite or) countable index set. We call 
$A=(a_{ij})_{i,j\in I}$ a \textit{Borcherds--Cartan\ matrix}
if the following three conditions are satisfied:
(1)\ $a_{ii}=2$ or $a_{ii}\in \mathbb{Z}_{\le 0}$ for each $i\in I$; 
(2)\ $a_{ij}\in \mathbb{Z}_{\le 0}$ for all $i,j\in I$ with $i\neq j$; 
(3)\ $a_{ij}=0$ if and only if $a_{ji}=0$ for all $i,j\in I$ with $i\neq j$.

An index $i\in I$ is said to be \textit{real} if $a_{ii}=2$, and 
\textit{imaginary} if $a_{ii}\le 0$. Denote by 
$I^{re}:=\{ i\in I\ |\ a_{ii}=2 \}$ the set of real indices, and by 
$I^{im}:=\{ i\in I\ |\ a_{ii}\le 0 \} =I\setminus I^{re}$ the set of 
imaginary indices. A Borcherds--Cartan matrix $A$ is said to be 
\textit{symmetrizable} if there exists a diagonal matrix 
$D=\mathrm{diag}(d_i ) _{i\in I}$, with $d_i \in \mathbb{Z}_{>0}$, 
such that $DA$ is symmetric. Also, if $a_{ii}\in 2\mathbb{Z}$ for all $i\in I$, then $A$ is said to be \textit{even}. 
Throughout this paper, we assume that the Borcherds--Cartan matrix is symmetrizable and even. 

For a given Borcherds--Cartan matrix $A=(a_{ij})_{i,j\in I}$, 
a \textit{Borcherds--Cartan datum} is a quintuple 
$\bigl( A,\varPi :=\{ \alpha _i \} _{i\in I}, \varPi ^{\vee}:=\{ \alpha _i ^{\vee} \} _{i\in I},
P, P^{\vee}\bigl)$, where $\varPi $ and $\varPi ^{\vee}$ are the sets of
\textit{simple roots} and \textit{simple coroots}, respectively, $P^{\vee}$ 
is a \textit{coweight lattice}, and $P:=\mathrm{Hom}_{\mathbb{Z}}(P^{\vee},\mathbb{Z})$ 
is a \textit{weight lattice}. We set $\mathfrak{h}:=P^{\vee}\otimes _{\mathbb{Z}}\mathbb{C}$,
and call it the \textit{Cartan subalgebra}. Let
$\mathfrak{h}^* :=\mathrm{Hom}_{\mathbb{C}}(\mathfrak{h},\mathbb{C})$ 
denote the full dual space of $\mathfrak{h}$.
In this paper, we assume that $\varPi ^{\vee} \subset \mathfrak{h}$
and $\varPi \subset \mathfrak{h}^*$ are both linearly independent over $\mathbb{C}$.
Let $P^+ :=\{ \lambda \in P \ |\ \alpha _i ^{\vee}(\lambda )
\ge 0 \ \mathrm{for\ all}\ i\in I \}$ be the set of dominant integral weights, and 
$Q:=\bigoplus _{i\in I}\mathbb{Z}\alpha _i $ the \textit{root lattice}; 
we set $Q^+ :=\sum _{i\in I}\mathbb{Z}_{\ge 0}\alpha _i $.
Let $\mathfrak{g}$ be the generalized Kac--Moody algebra associated with 
a Borcherds--Cartan datum $(A,\varPi ,\varPi ^{\vee},P,P^{\vee} )$. 
We have the root space decomposition
$\mathfrak{g}=\bigoplus _{\alpha \in \mathfrak{h}^*}
\mathfrak{g}_{\alpha }$, where $\mathfrak{g}_{\alpha }
:=\{ x\in \mathfrak{g}\ |\ \mathrm{ad}(h)(x)=h(\alpha )x \ 
\mathrm{for\ all}\ h\in \mathfrak{h} \},$ and $\mathfrak{h}=\mathfrak{g}_0.$
Denote by $\varDelta :=\{ \alpha \in \mathfrak{h}^* \mid 
\mathfrak{g}_{\alpha }\neq \{ 0 \} ,\ \alpha \neq 0 \}$ 
the set of roots, and by $\varDelta ^+$ the set of positive roots.
Note that $\varDelta \subset Q \subset P$ and 
$\varDelta =\varDelta ^+ \sqcup (- \varDelta ^+)$.
Also, we define a Coxeter group $\mathcal{W}_{re}:=\langle r_i \ |\ i\in I^{re} 
\rangle _{\mathrm{group}} \subset \mathrm{GL}(\mathfrak{h}^*)$,
where $r_i$ for $i\in I^{re}$ is a simple reflection,
and set $\varDelta _{im}:=\mathcal{W}_{re}\varPi _{im}$,
where $\varPi _{im}:=\{ \alpha _i \} _{i\in I^{im}}$ is the set of 
\textit{imaginary simple roots}; note that $\varDelta _{im} \subset \varDelta ^+$.
In addition, we set $\varDelta _{re} :=\mathcal{W}_{re}\varPi _{re}$ and 
$\varDelta _{re}^+ :=\varDelta _{re} \cap \varDelta ^+$, where 
$\varPi _{re}:=\{ \alpha _i \} _{i\in I^{re}}$ is the set of \textit{real simple roots}.
As in the case of ordinary Kac--Moody algebras, it is easily checked that the coroot 
$\beta ^{\vee}:=w\alpha _i ^{\vee}$ of $\beta =w\alpha _i \in \varDelta _{re}^+ \sqcup \varDelta _{im}$ 
is well-defined (see [\textbf{JL}, \S 2.1.9]).
\begin{define}\label{2.1.1}
Let $q$ be an indeterminate. The quantized universal enveloping algebra 
$U_q (\mathfrak{g})$ associated with a Borcherds--Cartan datum 
$(A,\varPi ,\varPi ^{\vee},P,P^{\vee} )$, with $D=\mathrm{diag}(d_i )_{i\in I}$ as above, 
is a $\mathbb{C}(q)$-algebra generated by the symbols 
$e_i $, $f_i$, $i\in I$, and $q^h$, $h\in P^{\vee}$, subject to the following relations:
\begin{itemize}
\item[$\cdot$]\ $q^0=1$, $q^{h_1}q^{h_2}=q^{h_1 +h_2}$ for $h_1 ,h_2 \in P^{\vee}$, 
\item[$\cdot$]\ $q^h e_i q^{-h} = q^{h(\alpha _i)} e_i$, 
$q^h f_i q^{-h} = q^{-h(\alpha _i)} f_i$ for $h\in P^{\vee}$ and $i\in I$, 
\item[$\cdot$]\ $[e_i ,f_j ]=\delta _{ij}
\displaystyle{\frac{K_i -K_i ^{-1}}{q_i -q_i ^{-1}}}$ for $i,j \in I$, where we set 
$q_i := q^{d_i }$ and $K_i := q^{d_i \alpha _i ^{\vee}}$, 
\item[$\cdot$]\ $\displaystyle{\sum _{r=0} ^{1-a_{ij}}} (-1)^r 
\begin{bmatrix}1-a_{ij} \\ r \end{bmatrix}_i e_i ^{1-a_{ij}-r} e_j e_i ^r =0$
for $i\in I^{re}$ and $j\in I$, with $i\neq j$, 
\item[$\cdot$]\ $\displaystyle{\sum _{r=0} ^{1-a_{ij}}} (-1)^r 
\begin{bmatrix}1-a_{ij} \\ r \end{bmatrix}_i f_i ^{1-a_{ij}-r} f_j f_i ^r =0$
for $i\in I^{re}$ and $j\in I$ , with $i\neq j$, 
\item[$\cdot$]\ $[e_i ,e_j ]=[f_i ,f_j ]=0$ if $a_{ij}=0$.
\end{itemize}
Here, we set $[n]_i :=\displaystyle{\frac{q_i ^n -q_i ^{-n}}{q_i -q_i ^{-1}}}$,
$[n]_i ! :=\displaystyle{\prod _{k=1}^n}[k]_i$, and 
$\begin{bmatrix}m \\ n \end{bmatrix}_i 
:=\displaystyle{\frac{[m]_i !}{[m-n]_i ! [n]_i !}}$ for $i\in I$.
\end{define}

If $a_{ii}<0$, then we set $c_i := -\frac{1}{2}a_{ii} \in 
\mathbb{Z}_{>0}$, and define 
$\{ n \} _i :=\frac{q_i ^{c_i n}-q_i ^{-c_i n}}
{q_i ^{c_i} -q_i ^{-c_i}}$. If $a_{ii} =0$, we set 
$\{ n \} _i :=n$. We define \textit{divided\ powers} by 
$e_i ^{(n)}:=\frac{e_i ^n}{[n]_i !}$, 
$f_i ^{(n)}:=\frac{f_i ^n}{[n]_i !}$ 
if $i\in I^{re}$, and by 
$e_i ^{(n)}:=e_i ^n$, $f_i ^{(n)}:= f_i ^n$ if $i\in I^{im}$.
Here we understand that $e_i ^{(0)}=f_i ^{(0)}:=1$, and 
$e_i ^{(n)}=f_i ^{(n)}:=0$ for $n<0$.
Let $U_q ^+ (\mathfrak{g})$ and $U_q ^- (\mathfrak{g})$
be the subalgebras of $U_q (\mathfrak{g})$ generated by 
$e_i $, $i\in I$, and $f_i $, $i\in I$, respectively.

Let $V(\lambda )$ be the irreducible highest weight 
$U_q (\mathfrak{g})$-module of highest weight $\lambda \in P^+$.
We define the \textit{Kashiwara\ operators} 
$\Tilde{e}_i , \Tilde{f}_i$, $i\in I$, on $V(\lambda )$ 
in the following way (see [\textbf{JKK}, \S 4]). 
Since we have $V(\lambda )=\bigoplus _{n \ge 0} f_i ^{(n)} \mathrm{Ker}(e_i)$ 
for each $i\in I$, and the weight space decomposition 
$V(\lambda ) = \bigoplus _{\mu \in P} V(\lambda )_{\mu}$, 
each weight vector $v \in V(\lambda )_{\mu }$, $\mu \in P$, 
has the unique expression $v = \sum _{n \ge 0} f_i ^{(n)} v_n $ such that 
(a) $v_n \in \mathrm{Ker}(e_i) \cap V(\lambda )_{\mu + n \alpha _i}$, 
(b) if $i\in I^{re}$ and $\alpha _i ^{\vee}(\mu + n\alpha _i) < n$, then $v_n =0$, and 
(c) if $i\in I^{im}$, $n>0$, and $\alpha _i ^{\vee}(\mu + n\alpha _i) = 0$, then $v_n =0$.
This expression for $v$ is called the \textit{$i$-string decomposition}.
Then, we define $\Tilde{f}_i v := \sum _{n \ge 0} f_i ^{(n+1)} v_n$
and $\Tilde{e}_i v := \sum _{n \ge 1} f_i ^{(n-1)} v_n$. 
Note that if $i\in I^{im}$, then $\Tilde{f}_i = f_i $.

Let $u_{\lambda} \in V(\lambda)_{\lambda } \setminus \{ 0 \}$ denote the highest weight vector, 
and $B(\lambda )$ the crystal basis of $V(\lambda )$ with the crystal lattice $L(\lambda )$. 
Here, $L(\lambda )$ is a free module over the local ring 
$\{ f(q)/g(q) \mid f(q), g(q) \in \mathbb{C}[q],\ g(0)\neq 0 \}$ generated by 
$\Tilde{f}_{i_l}\cdots \Tilde{f}_{i_1} u_{\lambda}$, 
$l\ge 0$, $i_1 , \ldots , i_l \in I$, and we have 
$B(\lambda) =\{ \Tilde{f}_{i_l}\cdots \Tilde{f}_{i_1} u_{\lambda}
\mod qL(\lambda) \mid l\ge 0,\ i_1 , \ldots , i_l \in I\} \setminus \{ 0\}$.
Let $\mathbf{A}:=\mathbb{C}[q,q^{-1}]$, and denote by 
$V(\lambda )^{\mathbf{A}}$ the $\mathbf{A}$-form of $V(\lambda )$.
Let $\bar{\ }:U_q (\mathfrak{g})\rightarrow U_q (\mathfrak{g})$ be the 
$\mathbb{C}$-algebra automorphism defined by 
$$q\longmapsto q^{-1},\ q^h\longmapsto q^{-h},\ 
e_i \longmapsto e_i ,\ f_i \longmapsto f_i,$$
for $h\in P^{\vee}$ and $i\in I$. 
Also, we define a $\mathbb{C}$-linear automorphism 
$\bar{\ }$ on $V(\lambda )$ by 
$\overline{X u_{\lambda}}:=\overline{X}u_{\lambda}$
for $X\in U_q (\mathfrak{g})$.
Let $\{ G_{\lambda}(b) \} _{b\in B(\lambda)}$ denote the 
global basis of $V(\lambda)$. 
We know from [\textbf{JKK}, Theorem 9.3] 
that the element $G_{\lambda}(b)$, $b\in B(\lambda)$, 
is characterized by the following three conditions: 
(i) $\overline{G_{\lambda}(b)}=G_{\lambda}(b)$, 
(ii) $G_{\lambda}(b) \in V(\lambda )^{\mathbf{A}} \cap L(\lambda)$, and 
(iii) $G_{\lambda}(b) \equiv b \mod qL(\lambda)$.

\subsection{The associated Kac--Moody algebra $\Tilde{\mathfrak{g}}$}
We associate a Cartan matrix $\widetilde{A}$ with 
a given Borcherds--Cartan matrix $A=(a_{ij})_{i,j\in I}$ as follows.
Set $\Tilde{I}:=\{ (i,1) \} _{i\in I^{re}} \sqcup 
\{ (i,m) \} _{i\in I^{im}, m\in \mathbb{Z}_{\ge 1}}$, and define a Cartan matrix 
$\widetilde{A}:=(\Tilde{a}_{(i,m), (j,n)})_{(i,m),(j,n) \in \Tilde{I}}$ by
\begin{center}
$\begin{cases}
\Tilde{a}_{(i,m), (i,m)} :=2 
\ \mathrm{for}\ (i,m) \in \Tilde{I}, \\
\Tilde{a}_{(i,m), (j,n)} :=a_{ij} \ 
\mathrm{for}\ (i,m),(j,n) \in \Tilde{I}\ 
\mathrm{with}\ (i,m) \neq (j,n).
\end{cases}$
\end{center}
Note that if $A$ is symmetrizable, then so is $\widetilde{A}$
(see [\textbf{I}, Lemma 4.1.1]).
Let us denote by $\bigl( \widetilde{A},\widetilde{\varPi }:=
\{ \Tilde{\alpha} _{(i,m)} \} _{(i,m) \in \Tilde{I}}, \widetilde{\varPi }^{\vee}:=
\{ \Tilde{\alpha} _{(i,m)} ^{\vee}\} _{(i,m) \in \Tilde{I}}, \widetilde{P}, 
\widetilde{P}^{\vee} \bigl)$ the Cartan datum associated with $\widetilde{A}$,
where $\widetilde{\varPi }, \widetilde{\varPi }^{\vee}$ are the sets of simple roots and simple coroots, respectively, 
$\widetilde{P}$ is a weight lattice, and $\widetilde{P}^{\vee}$ is a coweight lattice.
Let $\Tilde{\mathfrak{g}}$ be the associated Kac--Moody algebra, 
$\Tilde{\mathfrak{h}}:=\widetilde{P}^{\vee}\otimes _{\mathbb{Z}} \mathbb{C}$
the Cartan subalgebra, and $\widetilde{W}$ the Weyl group.
Note that every permutation on the subset $\{ (i,m) \} _{m\in \mathbb{Z}_{\ge 1}} \subset \Tilde{I}$ 
induces a (Dynkin) diagram automorphism of $\Tilde{\mathfrak{g}}$
for each $i\in I^{im}$. We denote by $\mathfrak{S}_i$ the permutation group 
on the subset $\{ (i,m) \} _{m\in \mathbb{Z}_{\ge 1}}\subset \Tilde{I}$ for each $i\in I^{im}$, 
and set $\Omega :=\prod _{i\in I^{im}} \mathfrak{S}_i$.

For notational simplicity, we write 
$(\mathbf{i,m}) := ( (i_s , m_s) ) _{s=1} ^k 
=( (i_k ,m_k ),\ldots , (i_2, m_2), (i_1, m_1) )$ if 
$\mathbf{i}=(i_k ,\ldots ,i_2 ,i_1) \in I^k$ and 
$\mathbf{m}=(m_k ,\ldots ,m_2 ,m_1) \in \mathbb{Z}^k$
for $k\ge 0$, and call it an \textit{ordered index} if 
$m_r =1$ for $i_r \in I^{re}$, and if $m_{x_s} =s$ for all 
$s=1,2,\ldots ,t$, where $\{ x_1 , x_2 , \ldots , x_t \} 
=\{ 1\le  x \le k \mid i_x = i \}$, with $1 \le x_1 < x_2 < \cdots < x_t \le k$, for each $i\in I^{im}$.
We set $\mathcal{I}:=\bigcup _{k=1} ^{\infty} I^k$, 
$\widetilde{\mathcal{I}}:=\bigcup _{k=1}^{\infty}
\Tilde{I}^k$, and denote by $\widetilde{\mathcal{I}}_{\mathrm{ord}}$ the set of 
all ordered indices. 
Since the $\mathbf{m}$ for which $\mathbf{(i,m)} \in \widetilde{\mathcal{I}}_{\mathrm{ord}}$
is determined uniquely by $\mathbf{i}$, we have a bijection 
$\mathcal{I}\rightarrow \widetilde{\mathcal{I}}_{\mathrm{ord}},\ \mathbf{i}\mapsto \mathbf{(i,m)}.$
Note that $\Omega $ acts on $\widetilde{\mathcal{I}}$ diagonally.
\begin{ex}\label{2.2.2}
If $I = \{ 1, \bm{2}, \bm{3} \}$, $I^{re} =\{ 1 \}$ and $I^{im} = \{ \bm{2}, \bm{3} \}$, 
then $\Tilde{I} = \{ (1,1) \} \sqcup \{ ( \bm{2}, m) ,$
$( \bm{3}, m) \} _{m\in \mathbb{Z}_{\ge 1}}$.
An ordered index $\mathbf{(i,m)} \in \widetilde{\mathcal{I}}_{\mathrm{ord}}$
corresponding to 
$\mathbf{i}=(\bm{3},1,\bm{3},\bm{3},1,\bm{2},\bm{3},1,\bm{3},\bm{2},1,\bm{2})\in \mathcal{I}$
is as follows:
\begin{center}
$\mathbf{(i,m)}=\bigl( (\bm{3},5),(1,1),(\bm{3},4),(\bm{3},3),(1,1),(\bm{2},3),(\bm{3},2),(1,1),
(\bm{3},1),(\bm{2},2),(1,1),(\bm{2},1) \bigl)$.
\end{center}
\end{ex}
We also define the subset 
$\widetilde{\mathcal{I}}_{\mathrm{gen}}\subset \widetilde{\mathcal{I}}$
of \textit{generic\ indices} in $\widetilde{\mathcal{I}}$ by
$\widetilde{\mathcal{I}}_{\mathrm{gen}}:=
\{ \mathbf{(i,m)}=((i_s , m_s)) _{s=1} ^k  \in \widetilde{\mathcal{I}} \mid 
m_s \neq m_t$ for $s\neq t$ with $i_s = i_t \in I^{im}\}$; 
note that $\widetilde{\mathcal{I}}_{\mathrm{ord}} \subset 
\widetilde{\mathcal{I}}_{\mathrm{gen}} \subset \widetilde{\mathcal{I}}$, 
and that $\widetilde{\mathcal{I}}_{\mathrm{gen}}=\Omega \ \! \widetilde{\mathcal{I}}_{\mathrm{ord}}$. 
In particular, $\widetilde{\mathcal{I}}_{\mathrm{gen}}$ is stable under the action of $\Omega$.

\subsection{The monoid $\mathcal{W}$}
In this subsection, we give a brief review of the monoid $\mathcal{W}$ 
and its properties. For more details, we refer the reader to [\textbf{I}, \S 2.2].

If we define $r_i \in \mathrm{GL}(\mathfrak{h}^*),\ i\in I^{im}$, by
$r_i (\mu ):=\mu -\alpha _i ^{\vee}(\mu )\alpha _i $ for $\mu \in \mathfrak{h}^*$,
then the inverse of $r_i $ is given by
$$r_i ^{-1}(\mu )=\mu +\frac{1}{1-a_{ii}}\alpha _i ^{\vee}(\mu )\alpha _i\ 
\mathrm{for}\ \mu \in \mathfrak{h}^*.$$
Note that this $r_i$ has an infinite order in $\mathrm{GL}(\mathfrak{h}^*)$.
\begin{define}\label{2.3.1}
$\mathrm{([\mathbf{I},\ Definition\ 2.2.1]).}$
Let $\mathcal{W}$ denote the monoid generated
by the symbols $\Tilde{r}_i ,\ i\in I$, 
subject to the following relations:
\begin{itemize}
\item[(1)]
$\Tilde{r}_i ^2 =1$ for all $i\in I^{re}$;
\item[(2)]
if $i,j \in I^{re},\ i\neq j,$ and the order of 
$r_i r_j \in \mathrm{GL}(\mathfrak{h}^*)$ is 
$m\in \{ 2,3,4,6 \}$, then we have 
\begin{center}
$(\Tilde{r}_i \Tilde{r}_j )^m 
=(\Tilde{r}_j \Tilde{r}_i )^m =1$;
\end{center}
\item[(3)]
if $i\in I^{im}$, then for all $j\in I\setminus \{ i\}$ such that 
$a_{ij}=0$, we have $\Tilde{r}_i \Tilde{r}_j =\Tilde{r}_j \Tilde{r}_i$.
\end{itemize}
\end{define}
Each element $w \in \mathcal{W}$ can be written as 
a product $w= \Tilde{r}_{i_1} \Tilde{r}_{i_2} \cdots \Tilde{r}_{i_k}$ 
of generators $\Tilde{r}_i$, $i \in I$. If the number $k$ is minimal among 
all the expressions for $w$ of the form above, then $k$ is called the 
\textit{length} of $w$ and the expression 
$\Tilde{r}_{i_1} \Tilde{r}_{i_2} \cdots \Tilde{r}_{i_k}$
is called a \textit{reduced\ expression}. In this case, we write $\ell (w) =k$.
Since the $r_i \in \mathrm{GL}(\mathfrak{h}^*)$, $i\in I$, satisfy the conditions
$(1),(2)$ and $(3)$ in Definition \ref{2.3.1}, we have the following 
(well-defined) homomorphism of monoids:
\begin{center}
$\mathcal{W}\longrightarrow \mathrm{GL}(\mathfrak{h}^*),
\ \Tilde{r}_i \longmapsto r_i , \  \mathrm{for} \ i\in I.$
\end{center}
For simplicity, we write $r_i $ for $\Tilde{r}_i$ in $\mathcal{W}$.
Remark that $\mathcal{W}_{re}$ and 
$\langle r_i \mid i\in I^{re} \rangle _{\mathrm{monoid}} \subset \mathcal{W}$
are isomorphic as groups, where $\langle r_i \mid i\in I^{re} \rangle _{\mathrm{monoid}}$
denotes the submonoid of $\mathcal{W}$ generated by $r_i $, $i\in I^{re}$; 
note that this submonoid is in fact a group by Definition \ref{2.3.1} (1). 
Hence we may (and do) regard the group
$\mathcal{W}_{re}$ as a submonoid of $\mathcal{W}$. 
For $\beta = w \alpha _i \in \varDelta _{re} ^+ \sqcup \varDelta _{im}$, 
$i\in I$, $w\in \mathcal{W}_{re}$, the element 
$r_{\beta} := w r_i w^{-1}$ is well-defined (see [\textbf{I}, \S 2.2]).

\begin{define}\label{2.3.2}
$(\mathrm{[\mathbf{I},\ Definition\ 2.2.6]).}$\ 
For $w\in \mathcal{W}$ and 
$\beta \in \varDelta _{re} ^+ \sqcup \varDelta _{im}$, we write $w \rightarrow r_{\beta} w$
if $\ell (r_{\beta } w)>\ell (w)$. Also, we define a partial order
$\le $ on $\mathcal{W}$ as follows: $w \le w' \ \mathrm{in}\ \mathcal{W}$ if there exist 
$w_0 ,w_1 ,\ldots ,w_l \in \mathcal{W}$
such that $w=w_0\rightarrow w_1 \rightarrow \cdots \rightarrow w_l =w' .$
\end{define}
Note that if $v \le w$ in $\mathcal{W}$ and 
$w=r_{i_l} \cdots r_{i_2} r_{i_1}$ is a reduced expression, 
then there exists a reduced expression $v= r_{i_{x_p}} \cdots r_{i_{x_2}} r_{i_{x_1}}$, 
$l \ge x_p > \cdots >x_2 > x_1 \ge 1$, 
by the \textit{Exchange\ Property} of $\mathcal{W}$ (see [\textbf{I}, \S 2.2]). 
However, the converse does not hold in general. For example, if $i, j \in I^{im}$ with $a_{ij} \neq 0$, 
then $r_i $ and $r_i r_j$ are not comparable even though $r_i$ is a subword of $r_i r_j$ and 
each of them is a (unique) reduced expression.

\begin{define}\label{2.3.3}
$(\mathrm{[\mathbf{I},\ Definition\ 2.2.10]).}$
Let $w=w_k r_{i_k} \cdots w_1 r_{i_1} w_0 \in \mathcal{W}$, with 
$i_1 ,\ldots ,i_k \in I^{im}$ and $w_0 ,w_1 ,\ldots ,w_k \in \mathcal{W}_{re}$, 
be a reduced expression, namely, $\ell (w) = \ell (w_k) + 1 + \cdots + \ell (w_1) + 1 + \ell (w_0)$.
We call this expression a dominant reduced expression if it satisfies 
\begin{center}
$r_{i_s} w_{s-1} r_{i_{s-1}} \cdots w_1 r_{i_1} w_0 (P^+) \subset P^+$
for all $s=1,2,\ldots ,k$. 
\end{center}
\end{define}
Note that every $w\in \mathcal{W}$ 
has at least one dominant reduced expression.
Indeed, if we choose an expression of the form above in such a way that
the sequence $\bigl( \ell (w_0 ), \ell (w_1), \ldots , \ell (w_k ) \bigl)$ 
is minimal in lexicographic order among all the reduced expressions 
of the form above, then it is a dominant reduced expression (see [\textbf{I}, \S 2.2]).

\subsection{Joseph--Lamprou's path model}
In this subsection, following [\textbf{JL}], we review Joseph--Lamprou's path model
for generalized Kac--Moody algebras. For more details, we refer the 
reader to [\textbf{JL}, \textbf{I}].

Let $\mathfrak{h}_{\mathbb{R}}$ denote a real form of $\mathfrak{h}$,
and $\mathfrak{h}_{\mathbb{R}} ^*$ its full dual space.
Let $\mathbb{P}$ be the set of all piecewise-linear continuous maps
$\pi :[0,1]\longrightarrow \mathfrak{h}_{\mathbb{R}} ^*$
such that $\pi (0)=0$ and $\pi (1) \in P$, 
where we set $[0,1]:=\{ t\in \mathbb{R}\ |\ 0\le t\le1 \}$. 
Also, we set $H_i ^{\pi }(t):=\alpha _i ^{\vee}\bigl( \pi (t) \bigl)$ for 
$t\in [0,1]$, and then $m_i ^{\pi }:=\min \{ H_i ^{\pi }(t) \mid H_i ^{\pi }(t)\in \mathbb{Z}, t\in [0,1]\}$.

Now we define the \textit{root operators} 
$e_i ,f_i :\mathbb{P}\longrightarrow \mathbb{P} \sqcup \{ \mathbf{0} \}$ for $i\in I$.
First, we set 
$f_+ ^i (\pi ):= \max \{ t\in [0,1] \mid H_i ^{\pi }(t)=m_i ^{\pi} \}$. If $f_+ ^i (\pi )<1$, then we can define 
$f_- ^i (\pi ) := \min \{ t\in [f_i ^+ (\pi ),1] \mid H_i ^{\pi }(t)=m_i ^{\pi }+1 \}$. In this case, we set 
\begin{align*}
(f_i \pi )(t):=
\begin{cases}
\pi (t) & t\in [0,f_+ ^i (\pi )], \\
\pi \bigl ( f_+ ^i (\pi ) \bigl ) + r_i \bigl ( \pi (t) -
\pi \bigl ( f_+ ^i (\pi ) \bigl ) \bigl ) & t\in [f_+ ^i (\pi), f_- ^i (\pi )] \\
\pi (t)-\alpha _i & t\in [f_- ^i (\pi ),1].
\end{cases}
\end{align*}
Otherwise (i.e., if $f_+ ^i (\pi )=1$), we set $f_i \pi :=\mathbf{0}$.

Next, we define the operator $e_i $ for $i\in I^{re}$.
Set $e _+ ^i (\pi ):= \min \{ t\in [0,1] \mid H_i ^{\pi }(t) 
=m_i ^{\pi } \}$. If $e_+ ^i (\pi )>0$, then we can define
$e_- ^i (\pi ):= \max \{ t\in [0,e_+ ^i (\pi )] \mid 
H_i ^{\pi }(t)= m_i ^{\pi } +1 \}$. In this case, we set
\begin{align*}
(e_i \pi )(t):=
\begin{cases}
\pi (t) & t\in [0,e_- ^i (\pi )], \\
\pi \bigl ( e_- ^i (\pi ) \bigl ) + r_i \bigl ( \pi (t) -
\pi \bigl ( e_- ^i (\pi ) \bigl ) \bigl )  & t\in [e_- ^i (\pi), e_+ ^i (\pi )] \\
\pi (t)+\alpha _i & t\in [e_+ ^i (\pi ),1].
\end{cases}
\end{align*}
Otherwise (i.e., if $e_+ ^i (\pi )=0$), we set $e_i \pi :=\mathbf{0}$.

Finally, we define the operator $e_i $ for $i\in I^{im}$.
Set $e_- ^i (\pi ):=f_+ ^i (\pi )$. If $e_- ^i (\pi )<1$
and there exists $t \in [e_- ^i (\pi ),1]$ such that 
$H_i ^{\pi }(t)\ge m_i ^{\pi }+1-a_{ii}$, then we can define
$e_+ ^i (\pi ):=\min \{ t\in [e_- ^i (\pi ),1] \mid 
H_i ^{\pi }(t)=m_i ^{\pi }+1 -a_{ii} \}$. We set $e_i \pi :=\mathbf{0}$ if $e_- ^i (\pi )=1$, or $e_- ^i (\pi )<1$ and 
$H_i ^{\pi }(t) < m_i ^{\pi } + 1 - a_{ii} $ for all 
$t\in [e_- ^i (\pi ),1]$, or $e_- ^i (\pi )<1$ and   
there exists $t\in [e_- ^i (\pi ),1]$ such that 
$H_i ^{\pi }(t)\ge m_i ^{\pi }+1-a_{ii}$ and 
$H_i ^{\pi }(s) \le m_i ^{\pi }-a_{ii}$ for some 
$s\in [e_+ ^i (\pi ),1]$. Otherwise, we set
\begin{align*}
(e_i \pi )(t):=
\begin{cases}
\pi (t) & t\in [0,e_- ^i (\pi )], \\
\pi \bigl ( e_- ^i (\pi ) \bigl ) + r_i ^{-1} \bigl ( \pi (t) -
\pi \bigl ( e_- ^i (\pi ) \bigl ) \bigl ) & t\in [e_- ^i (\pi), e_+ ^i (\pi )] \\
\pi (t)+\alpha _i & t\in [e_+ ^i (\pi ),1].
\end{cases}
\end{align*}

Let $\lambda \in P^+$. We write 
$\mu \ge \nu $ for $\mu ,\nu \in \mathcal{W}\lambda
:= \{ w\lambda \in P \mid w\in \mathcal{W} \} $
if there exists a sequence of elements 
$\mu =: \lambda _0 , \lambda _1 , \ldots , 
\lambda _{s-1}, \lambda_s := \nu$ in $\mathcal{W} \lambda $
and positive roots $\beta _1 ,\ldots ,\beta _s  
\in \varDelta _{re}^+ \sqcup \varDelta _{im}$
such that $\lambda _{i-1} = r_{\beta _i} \lambda _i$
and $\beta _i ^{\vee}(\lambda _i )>0$ for $i=1,\ldots ,s$.
This relation $\ge $ on $\mathcal{W} \lambda$ defines a partial order.
For $\mu , \nu \in \mathcal{W} \lambda $ and $\beta \in \varDelta _{re}^+ \sqcup \varDelta _{im}$, 
we write $\mu \xleftarrow{\beta } \nu$ if $\mu = r_{\beta } \nu$, $\beta ^{\vee}(\nu )>0 $,
and $\mu $ covers $\nu $ by this partial order.
Note that the direction of the arrow $\xleftarrow{\beta}$ defined above 
is opposite to that in [\textbf{JL}, \S 5.1.1]. 
\begin{define}\label{2.4.1}
$\mathrm{([\mathbf{JL},\ \S 5.2.1]}).$\ 
For a rational number $a\in (0,1]$ and 
$\mu ,\nu \in \mathcal{W}\lambda $ with $\mu \ge \nu $,
an $a$-chain for $(\mu ,\nu )$ is a sequence 
$\mu =:\nu _0 \xleftarrow{\beta _1 }\nu _1
\xleftarrow{\beta _2 } \cdots \xleftarrow{\beta _s }\nu _s :=\nu $
of elements in $\mathcal{W}\lambda $ such that for each $i=1,2,\ldots ,s,\ 
a\beta _i ^{\vee}(\nu _i)\in \mathbb{Z}_{>0}$
if $\beta _i \in \varDelta _{re} ^+,$ and 
$a\beta _i ^{\vee}(\nu _i)=1$ if $\beta _i \in \varDelta _{im}$.
\end{define}
\begin{define}\label{2.4.2}
$\mathrm{([\mathbf{JL},\ \S 5.2.2]).}$\ 
Let $\bm{\lambda }:=(\lambda _1 >\lambda _2 >\cdots >\lambda _s)$
be a sequence of elements in $\mathcal{W}\lambda $, and
$\bm{a}:=(0=a_0 <a_1 <\cdots <a_s =1)$ a sequence of rational numbers.
Then, the pair $\pi :=(\bm{\lambda }; \bm{a})$ is called a generalized 
Lakshmibai--Seshadri path (GLS path for short) of shape $\lambda $
if it satisfies the following conditions (called the chain condition):
$(\mathrm{i})$ there exists an $a_i$-chain for 
$(\lambda _i ,\lambda _{i+1})$ for each $i=1,2,\ldots ,s-1$;
$(\mathrm{ii})$ there exists a $1$-chain for $(\lambda _s ,\lambda ).$
\end{define}
In this paper, we think of the pair $\pi =(\bm{\lambda };\bm{a})$ 
as a path belonging to $\mathbb{P}$ by
$\pi (t):=\sum _{i=1}^{j-1}(a_i -a_{i-1})\lambda _i +(t-a_{j-1})\lambda _j$
for $a_{j-1}\le t\le a_j$ and $j=1,2,\ldots ,s.$
We denote by $\mathbb{B}(\lambda )$ the set of all GLS paths of shape $\lambda $.

Now, we define a crystal structure on $\mathbb{B}(\lambda )$. 
Let $\pi \in \mathbb{B}(\lambda )$. We set $\mathrm{wt}(\pi ):=\pi (1)\in P$.
For each $i\in I^{re}$, we set $\varepsilon _i (\pi ) := -m_i ^{\pi },\ 
\varphi _i (\pi ) := \alpha _i ^{\vee}\bigl( \mathrm{wt}(\pi )\bigl) -m_i ^{\pi }
=H_i ^{\pi }(1) -m_i ^{\pi }.$
Then, we have $\varepsilon _i (\pi )=\max \{ n\in \mathbb {Z}_{\ge 0} \ |\ e_i ^n \pi \in \mathbb{P} \}$, 
and $\varphi _i (\pi )=\max \{ n\in \mathbb {Z}_{\ge 0} \ |\ f_i ^n \pi \in \mathbb{P} \}$.
For each $i\in I^{im}$, we set $\varepsilon _i (\pi ):=0,\ 
\varphi _i (\pi ):= \alpha _i ^{\vee}\bigl( \mathrm{wt}(\pi )\bigl) .$
By the definitions, we have $\varphi _i (\pi )=\varepsilon _i (\pi )+
\alpha _i ^{\vee}\bigl( \mathrm{wt}(\pi )\bigl)$ for all $i\in I$.
Next, we define the \textit{Kashiwara operators} on $\mathbb{B}(\lambda )$.
We use the root operators $e_i ,\ i\in I^{re},$ and 
$f_i ,\ i\in I,$ on $\mathbb{P}$ as Kashiwara operators.
For $e_i ,\ i\in I^{im}$, we use the ``cutoff" of the
root operators $e_i ,\ i\in I^{im}$, on $\mathbb{P}\ \bigl( \supset \mathbb{B}(\lambda )\bigl)$,
that is, if $e_i \pi \notin \mathbb{B}(\lambda )$, then we set $e_i \pi :=\mathbf{0}$ in 
$\mathbb{B}(\lambda )$ even if $e_i \pi \neq \mathbf{0}$ in $\mathbb{P}$.
Thus, $\mathbb{B}(\lambda )$ is endowed with a crystal structure. 
From [\textbf{JL},\ Proposition\ 6.3.5], we have
$\mathbb{B}(\lambda )=\mathcal{F}\pi _{\lambda } \setminus \{ \mathbf{0} \}$,
where $\mathcal{F}$ is the monoid generated by the Kashiwara operators $f_i ,\ i\in I$.

\subsection{Embedding of path models} 
In this subsection, we give a brief review of the construction of an embedding of 
Joseph--Lamprou's path model for a given generalized Kac--Moody algebra $\mathfrak{g}$ 
into Littelmann's path model for an associated Kac--Moody algebra $\Tilde{\mathfrak{g}}$. 
For more details, we refer the reader to [\textbf{I}, \S 4]. 

Let $(A, \varPi , \varPi ^{\vee}, P, P^{\vee})$ be a Borcherds--Cartan datum, and 
$(\widetilde{A}, \widetilde{\varPi }, \widetilde{\varPi }^{\vee}, \widetilde{P}, \widetilde{P}^{\vee})$
be an associated Cartan datum as in \S 2.2. For each 
$\mu \in \mathcal{W}P^+ =\mathcal{W}_{re }P^+$, we take (and fix) an element
$\Tilde{\mu }\in \widetilde{P}$ such that
\begin{align}\label{eq1}
\Tilde{\alpha} _{(i,m)}^{\vee}(\Tilde{\mu })
=\alpha _i ^{\vee}(\mu )\ \mathrm{for\ all}\ (i,m)\in \Tilde{I}.
\end{align}
Let $\widetilde{\mathbb{B}}(\Tilde{\mu })$ be 
the set of all (G)LS paths of shape $\Tilde{\mu }$
for $\Tilde{\mathfrak{g}}$, and $\widetilde{\mathcal{F}}$
the monoid generated by the Kashiwara operators 
$f_{(i,m)}, (i,m) \in \Tilde{I}$. Following the notation of \S 2.2, 
we write $F_{\mathbf{i}}=f_{i_k }\cdots f_{i_2} f_{i_1} \in \mathcal{F}$ 
for $\mathbf{i}=( i_k ,\ldots ,i_2 ,i_1 ) \in \mathcal{I}$, and 
$F_{\mathbf{(i,m)}}=f_{(i_k , m_k )} \cdots f_{(i_2 , m_2)} f_{(i_1, m_1)}
\in \widetilde{\mathcal{F}}$ for $\mathbf{(i, m)}=((i_s , m_s)) _{s=1} ^k \in \widetilde{\mathcal{I}}$.
\begin{prop}\label{2.5.1}
$\mathrm{([\mathbf{I},\ Proposition\ 4.1.2]).}$\ 
For a dominant integral weight $\lambda \in P^+$, the map
\begin{center}
$\widetilde{\ \ }:\ \mathbb{B}(\lambda ) \longrightarrow \widetilde{\mathbb{B}}(\Tilde{\lambda }),\ 
\pi = F_{\mathbf{i}}\pi _{\lambda } \longmapsto \Tilde{\pi }:= 
F_{\mathbf{(i,m)}} \pi _{\Tilde{\lambda }},$
\end{center}
is well-defined and injective, where the $\mathbf{m}$ for which 
$\mathbf{(i,m)} \in \widetilde{\mathcal{I}}_{\mathrm{ord}}$ 
is determined uniquely by $\mathbf{i}\in \mathcal{I}$.
\end{prop}

Recall that $B(\lambda )$ denotes the crystal basis of 
the irreducible highest weight $U_q (\mathfrak{g})$-module $V(\lambda )$
of highest weight $\lambda \in P^+$. 
Since we know from [\textbf{I}, Theorem 6.1.1] 
that $\mathbb{B}(\lambda )\cong B(\lambda )$ as crystals, 
we obtain the following embedding by Proposition \ref{2.5.1}:
\begin{align}\label{eq2}
\widetilde{\ \ }:\  
B(\lambda ) \hookrightarrow \widetilde{B}(\Tilde{\lambda }),\ 
b=\widetilde{F}_{\mathbf{i}} u_{\lambda } \mapsto 
\Tilde{b}:= \widetilde{F}_{\mathbf{(i,m)}} \Tilde{u}_{\Tilde{\lambda }}, 
\ \text{for}\ \mathbf{(i,m)} \in \widetilde{\mathcal{I}}_{\mathrm{ord}},
\end{align}
where $\widetilde{B}(\Tilde{\lambda })$ denotes the crystal basis of 
the irreducible highest weight $U_q (\Tilde{\mathfrak{g}})$-module 
$\widetilde{V}(\Tilde{\lambda })$ of highest weight $\Tilde{\lambda }\in \widetilde{P}^+$,
$\widetilde{F}_{\mathbf{i}}:=\Tilde{f}_{i_k} \cdots \Tilde{f}_{i_2} \Tilde{f}_{i_1}$,
$\widetilde{F}_{\mathbf{(i,m)}}:=\Tilde{f}_{(i_k , m_k)} \cdots \Tilde{f}_{(i_2 , m_2)} \Tilde{f}_{(i_1 , m_1)}$
are monomials of the Kashiwara operators, and $u_{\lambda } \in V(\lambda )$, 
$\Tilde{u}_{\Tilde{\lambda }} \in \widetilde{V}(\Tilde{\lambda })$ are the highest weight vectors. 

\section{Some representation-theoretical results}
\subsection{Decomposition rules for $U_q (\mathfrak{g})$-modules 
in the category $\mathcal{O}_{\mathrm{int}}$}
In this subsection, we show the decomposition rules for $U_q (\mathfrak{g})$-modules 
stated in Theorems \ref{thm1} and \ref{thm2}. 
For this purpose, we recall the  following decomposition rules for path crystals.
\begin{thm}\label{3.1.1}
$\mathrm{([\mathbf{I},\ Theorem\ 7.1.3]).}$
Let $\lambda ,\mu \in P^+$. Then, we have an isomorphism of crystals:
$$\mathbb{B}(\lambda )\otimes \mathbb{B}(\mu )\ \cong \ 
\bigsqcup _
{\begin{subarray}{c}
\pi \in \mathbb{B}(\mu ) \\ 
\Tilde{\pi }\ \!\!:\ \!\! \Tilde{\lambda }
\text{-}\mathrm{dominant}
\end{subarray}}
\mathbb{B}\bigl( \lambda +\pi (1) \bigl).$$
Here, $\Tilde{\pi } \in \widetilde{\mathbb{B}}(\Tilde{\mu })$ denotes the 
image of $\pi \in \mathbb{B}(\mu )$ under the embedding 
$\mathbb{B}(\mu ) \hookrightarrow \widetilde{\mathbb{B}}(\Tilde{\mu })$,
and it is said to be $\Tilde{\lambda }$-dominant
if $\Tilde{\pi }(t)+\Tilde{\lambda }$ belongs to the dominant Weyl chamber 
of $\Tilde{\mathfrak{g}}$ for all $t\in [0,1]$.
\end{thm}

Let $S \subset I$ be a subset. We set
$S^{re}:=S\cap I^{re}$ and $S^{im}:=S\cap I^{im}$.
Also, we set $\widetilde{S}:=\{ (i,1) \} _{i\in S^{re}}\sqcup \{ (i,m) \} _{i\in S^{im}, m\in \mathbb{Z}_{\ge 1}}$.
Let us denote by $\mathfrak{g}_S$ (resp., $\Tilde{\mathfrak{g}}_{\widetilde{S}}$)
the Levi subalgebra of ${\mathfrak{g}}$ corresponding to $S$ 
(resp., the Levi subalgebra of $\mathfrak{\Tilde{g}}$ corresponding to $\widetilde{S}$), 
and denote by $\mathbb{B}_S (\lambda )$ the set of all GLS paths of shape $\lambda $ for $\mathfrak{g}_S$.
\begin{thm}\label{3.1.2}
$\mathrm{([\mathbf{I},\ Theorem\ 7.2.2]).}$
Let $\lambda \in P^+$. Then, we have an isomorphism of $\mathfrak{g}_S$-crystals:
$$\mathbb{B}(\lambda )\ \cong \ \bigsqcup _
{\begin{subarray}{c}\pi \in \mathbb{B}(\lambda) \\ 
\Tilde{\pi }\ \!\! :\ \!\! \Tilde{\mathfrak{g}}_{\widetilde{S}}
\text{-}\mathrm{dominant}
\end{subarray}} \mathbb{B}_S \bigl( \pi (1) \bigl) .$$
Here, $\Tilde{\pi } \in \widetilde{\mathbb{B}}(\Tilde{\lambda })$ 
denotes the image of $\pi \in \mathbb{B}(\lambda )$ under the 
embedding $\mathbb{B}(\lambda ) \hookrightarrow \widetilde{\mathbb{B}}(\Tilde{\lambda })$,
and it is said to be $\Tilde{\mathfrak{g}}_{\widetilde{S}}$-dominant
if $\Tilde{\pi }(t)$ belongs to the dominant Weyl chamber of $\Tilde{\mathfrak{g}}_{\widetilde{S}}$
for all $t\in [0,1]$.
\end{thm}

From [\textbf{JKK}, Theorems 3.7 and 7.1], we know the existence and uniqueness of the crystal basis of 
$U_q (\mathfrak{g})$-modules in the category $\mathcal{O}_{\mathrm{int}}$, 
and the complete reducibility for $U_q (\mathfrak{g})$-modules in the category $\mathcal{O}_{\mathrm{int}}$
(see [\textbf{JKK}, Definition 3.1] for the definition of the category $\mathcal{O}_{\mathrm{int}}$).
Since $V(\lambda) \otimes V(\mu ),\ \lambda , \mu \in P^+$, belongs to 
$\mathcal{O}_{\mathrm{int}}$ for $\mathfrak{g}$, and $V(\lambda)$ belongs to 
$\mathcal{O}_{\mathrm{int}}$ for $\mathfrak{g}_S$ (as a $U_q (\mathfrak{g}_S)$-module),
Theorems \ref{thm1} and \ref{thm2} follow immediately 
from Theorems \ref{3.1.1} and \ref{3.1.2}.

\subsection{Analog of the Parthasarathy--Ranga Rao--Varadarajan conjecture for 
generalized Kac--Moody algebras}

In this subsection, we prove an analog of the Parthasarathy--Ranga Rao--Varadarajan 
conjecture for generalized Kac--Moody algebras.

In what follows, we denote by $[\mu ]$ the unique element in 
$\mathcal{W}_{re} \mu \cap P^+$ for $\mu \in \mathcal{W}P^+$.
By an argument similar to the one for [\textbf{Li1}, Proposition 7.1], we can 
show the following lemma.
\begin{lem}\label{3.2.1}
Let $\lambda $, $\mu \in P^+$, and 
$\pi :=(\mu _1 ,\ldots ,\mu _l ;\ \! a_0 ,a_1 ,\ldots ,a_l)\in \mathbb{B}(\mu )$.
If $\lambda +\pi (a_p)\in P^+$ for all $0\le p\le l-1$,
then there exists a $\lambda $-dominant path $\pi ' \in \mathbb{B}(\mu )$
such that $\lambda +\pi '(1)=[\lambda +\pi (1)]$.
\end{lem}

Let $r_{(i,m)}$, $(i,m) \in \Tilde{I}$, denote the simple reflection 
of the Weyl group $\widetilde{W}$ of $\Tilde{\mathfrak{g}}$, 
$\widetilde{W}_{re}\subset \widetilde{W}$ the subgroup of 
$\widetilde{W}$ generated by  $r_{(i,1)}$, $i\in I^{re}$, and 
$\widetilde{W}_{\mathrm{gen}}:=\{ R_{\mathbf{(i,m)}} \in \widetilde{W} \mid
\mathbf{(i,m)} \in \widetilde{\mathcal{I}}_{\mathrm{gen}} \}$,  where 
$R_{\mathbf{(i,m)}} := r_{(i_l , m_l)} \cdots r_{(i_1 , m_1)}$ for 
$\mathbf{(i,m)}=(( i_s , m_s))_{s=1} ^l \in \widetilde{\mathcal{I}}$.
Note that $\widetilde{W}_{\mathrm{gen}}$ is not closed under multiplication, 
and that $\widetilde{W}_{re}$ is isomorphic as a group to the submonoid 
$\mathcal{W}_{re}$ of $\mathcal{W}$; we write this isomorphism as 
\begin{align}\label{eq3}
\widetilde{\ \ }:\ \mathcal{W}_{re}\xrightarrow{\ \cong \ } \widetilde{W}_{re},\ 
w = r_{i_l} \cdots r_{i_1} \longmapsto 
\Tilde{w} := r_{(i_l ,1)} \cdots r_{(i_1 ,1)},\ 
\mathrm{for}\ i_1 , \ldots , i_l \in I^{re}.
\end{align}
Also, we can easily check that the map 
$\widetilde{W}_{\mathrm{gen}} \rightarrow \mathrm{GL(\mathfrak{h}^*)},
\ R_{\mathbf{(i,m)}} \mapsto R_{\mathbf{i}}$, is well-defined, 
where $R_{\mathbf{i}} := r_{i_l } \cdots r_{i_1 }$ for $\mathbf{i}=( i_l ,\ldots , i_1) \in \mathcal{I}$.
\begin{lem}\label{3.2.2}
Let $\lambda \in P^+$. If $F_{\mathbf{(i,m)}}\pi _{\Tilde{\lambda }}$, 
$\mathbf{(i,m)}\in \widetilde{\mathcal{I}}_{\mathrm{gen}}$,  
is not $\mathbf{0}$ in $\widetilde{\mathbb{B}}(\Tilde{\lambda })$,
then $F_{\mathbf{(i,\bar{m})}} \pi _{\Tilde{\lambda }}$ 
is not $\mathbf{0}$ in $\widetilde{\mathbb{B}}(\Tilde{\lambda })$,
where the element $\mathbf{(i,\bar{m})} \in \widetilde{\mathcal{I}}_{\mathrm{ord}}$
is determined uniquely by $\mathbf{i}\in \mathcal{I}$.
\end{lem}
\textbf{Proof.}\ 
Since $\mathbf{(i, m)}\in \widetilde{\mathcal{I}}_{\mathrm{gen}}$, 
there exists a permutation $\omega \in \Omega$ on $\widetilde{\mathcal{I}}$ 
that sends $\mathbf{(i,m)}$ to $\mathbf{(i,\bar{m})}$ (see \S 2.2). 
Hence we have a permutation on the path crystal $\widetilde{\mathbb{B}}(\Tilde{\lambda })$, 
induced by the diagram automorphism of $\Tilde{\mathfrak{g}}$ corresponding to the
$\omega \in \Omega$, that sends $F_{\mathbf{(i,m)}}\pi _{\Tilde{\lambda }}$ 
to $F_{\mathbf{(i,\bar{m})}}\pi _{\Tilde{\lambda }}$ (see [\textbf{NS}, Lemma 3.1.1]). 
Therefore, we see that $F_{\mathbf{(i,\bar{m})}}\pi _{\Tilde{\lambda }} \neq \mathbf{0}$. \qed

\vspace{2mm}
For a path $\eta = F_{\mathbf{(i, m)}} \pi _{\Tilde{\lambda }}
\in \widetilde{\mathbb{B}}(\Tilde{\lambda })$ with $\mathbf{(i,m)} \in \widetilde{\mathcal{I}}_{\mathrm{gen}}$, 
we set $\bar{\eta }:=F_{\mathbf{(i,\bar{m})}} \pi _{\Tilde{\lambda }}$, with 
$\mathbf{(i,\bar{m})}\in \widetilde{\mathcal{I}}_{\mathrm{ord}}$. 
From the proof of Lemma \ref{3.2.2}, $\bar{\eta}$ is independent 
of the choice of an expression for $\eta $ of the form 
$F_{\mathbf{(i, m)}} \pi _{\Tilde{\lambda }}$ with $\mathbf{(i,m)} \in \widetilde{\mathcal{I}}_{\mathrm{gen}}$.
The following lemma is obvious.

\begin{lem}\label{3.2.3}
Let $\eta \in \widetilde{\mathbb{B}}(\Tilde{\lambda })$, 
and write $\eta (1)$ as $\Tilde{\lambda }-\sum c_{(i,m)} \Tilde{\alpha}_{(i,m)},
c_{(i,m)} \in \mathbb{Z}_{\ge 0}$. If $c_{(i,m)} \in \{ 0,1 \}$ for all 
$i\in I^{im}$ and $m\in \mathbb{Z}_{\ge 1}$, then there exists a path 
$\pi \in \mathbb{B}(\lambda )$ such that $\Tilde{\pi }=\bar{\eta }$
in $\widetilde{\mathbb{B}}(\Tilde{\lambda })$.
\end{lem}

For $\mathbf{i}=(i_k ,i_{k-1},\ldots ,i_1) \in \mathcal{I}$,
we set $\mathbf{i}_{[s]} :=(i_s ,i_{s-1},\ldots ,i_1),\ 1\le s \le k$.
Let $\lambda ,\mu \in P^+$ and $R_{\mathbf{i}}, R_{\mathbf{j}}\in \mathcal{W},\ 
\mathbf{i}=(i_k ,\ldots ,i_2 ,i_1)$, $\mathbf{j}=(j_l ,\ldots ,j_2 ,j_1)
\in \mathcal{I}$, be such that $\alpha _{i_s}^{\vee}(R_{\mathbf{i}_{[s-1]}}\lambda )=1\ 
\mathrm{\bigl( resp.,}\ \alpha _{j_t}^{\vee}(R_{\mathbf{j}_{[t-1]}} \mu )=1 \bigl)$
if $i_s \in I^{im}\ \mathrm{(resp.,}\ j_t \in I^{im})$.
We take sequences $\mathbf{m}=(m_k ,\ldots ,m_2 ,m_1)$
and $\mathbf{n}=(n_l ,\ldots ,n_2 ,n_1)$ of positive integers such that 
$\mathbf{(i,m)}, \mathbf{(j,n)} \in \widetilde{\mathcal{I}}_{\mathrm{gen}}$,
and such that $m_s \neq n_t$ if $i_s =j_t$ in $I^{im}$.
Also, we set $\nu :=R_{\mathbf{i}}\lambda + R_{\mathbf{j}}\mu \in P$ 
and $\bar{\nu }:= R_{\mathbf{(i,m)}} \Tilde{\lambda }+
R_{\mathbf{(j,n)}} \Tilde{\mu }\in \widetilde{P}$;
note that $\bar{\nu }\neq \Tilde{\nu }$ in general, where $\Tilde{\nu } \in \widetilde{P}$
is defined by equation (\ref{eq1}) in \S 2.5 for a weight $\nu \in \mathcal{W}P^+$.
With this notation, we have the following.
\begin{thm}\label{3.2.4}
If $\bar{\nu } \in \widetilde{P}^+$, then 
$V(\nu )$ appears in the irreducible decomposition of $V(\lambda )\otimes V(\mu )$.
\end{thm}
\textbf{Proof.}
If we take $M \in \mathbb{Z}_{\ge 1}$ such that 
$(i,M)$ does not appear in $\mathbf{(i,m)}$ and $\mathbf{(j,n)}$ for any $i\in I$, 
then the equality $\alpha _i ^{\vee}(\nu )=\Tilde{\alpha}_{(i,M)}^{\vee}(\bar{\nu })$ 
holds for all $i \in  I$ by the definition of $\bar{\nu}$. 
Since $\Tilde{\alpha}_{(i,M)}^{\vee}(\bar{\nu }) \ge 0$, 
we deduce that $\nu $ is a dominant integral weight.
By Theorem \ref{thm1}, it suffices to show that
there exists $\pi \in \mathbb{B}(\mu )$ such that 
$\Tilde{\pi }\in \widetilde{\mathbb{B}}(\Tilde{\mu })$
is $\Tilde{\lambda }$-dominant and $\lambda +\pi (1)=\nu $.
Set $\eta _1 :=\bigl( R_{\mathbf{(i,m)}}^{-1} R_{\mathbf{(j,n)}} \Tilde{\mu };\ \! 0,1\bigl) \in 
\widetilde{\mathbb{B}}(\Tilde{\mu })$, 
with $R_{\mathbf{(i,m)}}, R_{\mathbf{(j,n)}} \in \widetilde{W}_{\mathrm{gen}}$.
By Lemma \ref{3.2.1}, there exists a $\Tilde{\lambda }$-dominant path
$\eta _2 \in \widetilde{\mathbb{B}}(\Tilde{\mu })$ such that 
$\Tilde{\lambda }+\eta _2 (1)=[\Tilde{\lambda }+\eta _1 (1)]=\bar{\nu }$.
If we write $\eta _2 (1)$ as $\Tilde{\mu }-\sum c_{(i,m)} \Tilde{\alpha} _{(i,m)}$,
then $c_{(i,m)} \in \{ 0,1\}$ for all $i\in I^{im}$ and $m\in \mathbb{Z}_{\ge 1}$
by the definition of $\bar{\nu}$.
Therefore, by Lemma \ref{3.2.3}, there exists 
$\pi \in \mathbb{B}(\mu )$ such that $\Tilde{\pi }=\bar{\eta }_2$.
From this, we conclude that $\Tilde{\pi }\in \widetilde{\mathbb{B}}(\Tilde{\mu })$
is $\Tilde{\lambda }$-dominant and $\lambda +\pi (1)=\nu $. \qed
\begin{cor}\label{3.2.5}
Let $\lambda ,\mu \in P^+$, $w_1 ,w_2 \in \mathcal{W}_{re}$, 
and set $\nu =w_1 \lambda +w_2 \mu $. 
If $\nu \in P^+$, then $V(\nu )$ appears in the irreducible decomposition of 
$V(\lambda )\otimes V(\mu )$.
\end{cor}
\begin{rem}\label{3.2.6}
For general elements $w_1, w_2$ in the monoid $\mathcal{W}$,
the statement of the corollary above is false.
To see this, we take $\lambda ,\mu \in P^+$ and 
$i\in I^{im}$ such that $\alpha _i ^{\vee}(\lambda )=0$ and 
$\alpha _i ^{\vee}(\mu )=1$. 
Then, $\nu :=\lambda +r_i \mu =\lambda +\mu -\alpha _i$
is a dominant integral weight for $\mathfrak{g}$. 
However, we have $\bigl( \mathbb{B}(\lambda )\otimes \mathbb{B}(\mu ) \bigl) _{\nu }
=\bigl \{ \pi _{\lambda }\otimes f_i \pi _{\mu } \bigl \}$, 
and $e_i (\pi _{\lambda }\otimes f_i \pi _{\mu })
=\pi _{\lambda }\otimes e_i f_i \pi _{\mu }
=\pi _{\lambda }\otimes \pi _{\mu } \neq \mathbf{0}$.
Therefore, there is no highest weight vector of weight $\nu $,
and $V(\nu )$ cannot be an irreducible component 
of $V(\lambda )\otimes V(\mu )$. 
As for $\Tilde{\mathfrak{g}}$, we have 
$\bar{\nu }:=\Tilde{\lambda }+r_{(i,1)} \Tilde{\mu }$, 
$\Tilde{\alpha} _{(i,1)}^{\vee}(\bar{\nu })
=\Tilde{\alpha} _{(i,1)}^{\vee}(\Tilde{\lambda })+
\Tilde{\alpha} _{(i,1)}^{\vee}(r_{(i,1)} \Tilde{\mu })
=\alpha _i ^{\vee}(\lambda )-\alpha _i ^{\vee}(\mu )=-1$, 
and hence we see that $\bar{\nu }$ is not a 
dominant integral weight for $\Tilde{\mathfrak{g}}$.
\end{rem}

\section{Demazure modules for generalized Kac--Moody algebras}
\subsection{Definition of Demazure modules}
Let $\lambda \in P^+$ and $w=v_l r_{j_l} \cdots v_1 r_{j_1} v_0 \in \mathcal{W}$, with 
$j_1 ,\ldots ,j_l \in I^{im},\ v_0 ,\ldots ,v_l \in \mathcal{W}_{re}$, be such that
\begin{align}\label{eq4}
\alpha _{j_s} ^{\vee}\bigl( v_{s-1} r_{j_{s-1}} \cdots 
v_1 r_{j_1} v_0 (\lambda ) \bigl) =1\ \mathrm{for\ all}\ s=1,2,\ldots ,l. 
\end{align}
Note that the condition $(\ref{eq4})$ for $w$ and $\lambda $ 
is independent of the choice of an expression for $w$
since this condition is preserved under the relations 
$(1), (2)$, and $(3)$ of Definition \ref{2.3.1}.
For such an element $w\in \mathcal{W}$, we set 
$V_w (\lambda ):=U_q ^+ (\mathfrak{g}) \ \! V(\lambda )_{w\lambda }$, 
where $V(\lambda )_{\mu } \subset V(\lambda)$ 
denotes the weight space of weight $\mu \in P$.
We call $V_w (\lambda )$ the \textit{Demazure\ (sub)module}
of $V(\lambda )$ of lowest weight $w\lambda $. 
Since $V_w (\lambda )$ depends only on the weight $w\lambda \in P$, 
we may assume that the element $w=R_{\mathbf{i}}$, 
with $\mathbf{i}=(i_k ,\ldots ,i_2 ,i_1)$, satisfies the following condition
(recall the notation $\mathbf{i}_{[t]} =(i_t ,\ldots ,i_2 ,i_1),\ t=1,\ldots ,k$):
\begin{align}\label{eq5}
\alpha _{i_t } ^{\vee}(R_{\mathbf{i}_{[t-1]}}\lambda )>0
\ \mathrm{for\ all\ }t=1,2,\ldots ,k.
\end{align}

We will show that the Demazure module $V_w (\lambda)$ is generated 
by a single weight vector. Namely, we will show that $\dim V(\lambda)_{w\lambda}=1$;
note that the action of $\mathcal{W}$ on the set of weights of 
$V(\lambda )$ does not necessarily preserve the weight multiplicities.
Before doing this, we recall the action of $\mathcal{W}_{re}$ on $\mathbb{B}(\lambda)$
(see [\textbf{JL}, \S 9.2]): for $\pi \in \mathbb{B}(\lambda )$ and 
$v=r_{i_k}\cdots r_{i_2} r_{i_1}\in \mathcal{W}_{re}$, with 
$i_1 , i_2 , \ldots ,i_k \in I^{re}$, define 
\begin{align}\label{eq6}
S_v \pi :=x_{i_k }^{\alpha _{i_k} ^{\vee}(R_{\mathbf{i}_{[k-1]}}\mathrm{wt}(\pi ))}
\cdots x_{i_2} ^{\alpha _{i_2}^{\vee}(R_{\mathbf{i}_{[1]}}\mathrm{wt}(\pi ))}
x_{i_1}^{\alpha _{i_1} ^{\vee}(\mathrm{wt}(\pi ))}\pi ,
\end{align}
where $x_i ^a :=f_i ^a$ if $a\ge 0,\ x_i ^a :=e_i ^{-a}$ if $a<0$, 
and  $\mathbf{i}=(i_k ,\ldots ,i_2 ,i_1)$.
It is obvious that $\mathrm{wt}(S_v \pi )=v\bigl( \mathrm{wt}(\pi )\bigl)$.
\begin{lem}\label{4.1.1}
If $w\in \mathcal{W}$ and $\lambda \in P^+$ satisfy condition (\ref{eq4}), then we have
$\dim V(\lambda )_{w\lambda }=1$.
\end{lem}
\textbf{Proof.}
Let $w=v_l r_{j_l} \cdots v_1 r_{j_1} v_0$, with $v_0 ,\ldots ,v_ l \in \mathcal{W}_{re}$, 
$j_1 ,\ldots ,j_l \in I^{im}$, be a reduced expression satisfying condition (\ref{eq4}). 
It suffices to show that $\# \mathbb{B}(\lambda )_{w\lambda }=1$.
As remarked above, we may assume further that this expression 
for $w$ satisfies condition (\ref{eq5}).

If we set $\pi := S_{v_l} f_{j_l} \cdots S_{v_1} f_{j_1}S_{v_0}\pi _{\lambda }$, 
then the $x$'s in the expression (\ref{eq6}) for $S_{v_s},\ s=0,1,\ldots ,l,$
are all $f$ by condition (\ref{eq5}), and hence $\pi =\pi _{w\lambda }$ 
by condition (\ref{eq4}). Also, by using the fact [\textbf{I}, Lemma 4.1.5], we can deduce 
that the element $\Tilde{\pi } \in \widetilde{\mathbb{B}}(\Tilde{\lambda })$
(see Proposition \ref{2.5.1}) can be written as 
$\Tilde{\pi }=S_{\Tilde{v}_l} f_{(j_l , m_l)} \cdots S_{\Tilde{v}_1} f_{(j_1 , m_1)}
S_{\Tilde{v}_0}\pi _{\Tilde{\lambda }}=\pi _{\Tilde{w}\Tilde{\lambda }}$, where 
$((j_s , m_s))_{s=1} ^l \in \widetilde{\mathcal{I}} _{\mathrm{ord}}$, 
$\Tilde{w}:=\Tilde{v}_l r_{(j_l , m_l)}  \cdots \Tilde{v}_1 r_{(j_1 , m_1)} \Tilde{v}_0
\in \widetilde{W}$, and each $\Tilde{v}_s \in \widetilde{W}_{re}$ corresponds to 
$v_s \in \mathcal{W}_{re}$ via the isomorphism (\ref{eq3}) in \S 3.2.
This shows that the image of $\mathbb{B}(\lambda )_{w\lambda }$
under the embedding of Proposition \ref{2.5.1} is contained in 
$\widetilde{\mathbb{B}}(\Tilde{\lambda })_{\Tilde{w}\Tilde{\lambda }}$; 
note that if $\pi _1$ and $\pi _2$ in $\mathbb{B}(\lambda )$
satisfy $\mathrm{wt}(\pi _1) = \mathrm{wt}(\pi _2)$, then 
$\mathrm{wt}( \Tilde{\pi} _1) = \mathrm{wt}( \Tilde{\pi} _2)$.
Since $\widetilde{\mathbb{B}}(\Tilde{\lambda })_{\Tilde{w}\Tilde{\lambda }}
=\{ \pi _{\Tilde{w} \Tilde{\lambda}} \}$, we conclude that 
$\mathbb{B}(\lambda )_{w\lambda }=\{ \pi _{w\lambda} \}$, 
and hence $\# \mathbb{B}(\lambda )_{w\lambda } =1$. \qed
\begin{rem}\label{4.1.2}
Since the element $\Tilde{w} \in \widetilde{W}$ 
in the proof of Lemma \ref{4.1.1} is a minimal coset representative 
of a coset in $\widetilde{W}/ \widetilde{W}_{\Tilde{\lambda}}$, 
where $\widetilde{W}_{\Tilde{\lambda }} \subset \widetilde{W}$ denotes the 
isotropy subgroup for $\Tilde{\lambda } \in \widetilde{P}^+$
(see $[\mathbf{Hu},\ \S 1.10 \ \mathrm{and}\ \S 1.12]$), 
by the uniqueness of minimal coset representative,
$\Tilde{w}$ is independent of the choice of an expression for $w$ 
satisfying conditions (\ref{eq4}) and (\ref{eq5}) for $\lambda \in P^+$. 
\end{rem}

\subsection{Minimal dominant reduced expressions}
In this subsection, we introduce a specific expression for $w \in \mathcal{W}$,
which satisfies conditions (\ref{eq4}) and (\ref{eq5}) for $\lambda \in P^+$, 
in order to state Theorem \ref{thm4} (see \S 1).
\begin{lem}\label{4.2.1}
Let $w \in \mathcal{W}$ satisfy conditions (\ref{eq4}) and (\ref{eq5}) for $\lambda \in P^+$. 
Then, $w$ has an expression $w=w_k r_{i_k}^{a_k} \cdots w_1 r_{i_1}^{a_1} w_0$, 
with $w_0 ,\ldots ,w_k \in \mathcal{W}_{re}$ and
$a_1 ,\ldots ,a_k \in \mathbb{Z}_{\ge 1}$, where $i_1 ,\ldots ,i_k \in I^{im}$ are all distinct.
Moreover, if $s\in \{ 1\le t \le k \mid a_t >1 \}$, then $a_{i_s , i_s}=0$, i.e., 
the $(i_s , i_s)$-entry of the Borcherds--Cartan matrix $A$ is zero.
\end{lem}
\textbf{Proof.}\ 
Let $w=v_l r_{j_l} \cdots v_1 r_{j_1} v_0$, with 
$v_0 , \ldots ,v_l \in \mathcal{W}_{re}$, 
$j_1 ,\ldots ,j_l \in I^{im}$, be a reduced expression 
satisfying conditions (\ref{eq4}) and (\ref{eq5}) for $\lambda \in P^+$.
If $j_s =j_t$ with $s>t$, then we have
$$\alpha _{j_s}^{\vee}\bigl( v_{s-1} r_{j_{s-1}} \cdots 
v_t r_{j_t} (v_{t-1} r_{j_{t-1}} \cdots v_1 r_{j_1} v_0 (\lambda) ) \bigl) 
=\alpha _{j_t}^{\vee} \bigl( v_{t-1} r_{j_{t-1}} \cdots v_1 r_{j_1} v_0 (\lambda ) \bigl) $$
by condition (\ref{eq4}). Also, by condition (\ref{eq5}), $r_{j_s} $ commutes with 
$v_{s-1}, \ldots , v_{t+1}, v_t$ and $r_{j_{s-1}} , \ldots , r_{j_{t+1}}, r_{j_t}$ since 
$\alpha _{j_s} ^{\vee}$ is an anti-dominant integral coweight by [\textbf{JL}, Lemma 2.1.11]. 
In particular, we obtain $a_{j_s , j_s}=0$. \qed 

\vspace{2mm}
We fix an expression $w=w_k r_{i_k}^{a_k} \cdots w_1 r_{i_1 }^{a_1} w_0$
given in Lemma \ref{4.2.1} for which the sequence
$\bigl( \ell(w_0), \ell(w_1), \ldots , \ell(w_k) \bigl)$
is minimal in lexicographic order among all such expressions of $w$. 
Then, this is a dominant reduced expression (see Definition \ref{2.3.3}).
We call this expression a \textit{minimal\ dominant\ reduced\ expression}
(with respect to $\lambda \in P^+$).

Here we collect some fundamental properties of minimal dominant reduced expression
$w=w_k r_{i_k}^{a_k} \cdots w_1 r_{i_1 }^{a_1} w_0$; 
these properties follow directly from the definition and by induction on $k$ and $\ell (w)$.
For each $0\le s \le k$, let 
$w_s = r_{s, \ell _s} \cdots r_{s, 2} r_{s, 1}$, with 
$r_{s, p} := r_{\alpha _{s, p}}$, $\alpha _{s, p} \in \varPi _{re}$, 
be a (fixed) reduced expression, where $\ell _s := \ell (w_s)$.
\begin{lem}\label{4.2.2}
With the notation above, the following statements hold.
\begin{itemize}
\item[(a)]
$\alpha _{s,p}^{\vee}(r_{s,p-1} \cdots r_{s, 2} r_{s, 1} r_{i_{s-1}} ^{a_{s-1}}
\cdots w_1 r_{i_1 }^{a_1} w_0 \lambda ) =1$
for all $0\le s \le k-1$ and $1\le p\le \ell _s$. 
\item[(b)]
$\alpha _{s,p}^{\vee}(\alpha _{s,p-1}) =-1$
for all $0\le s \le k-1$ and $1< p \le \ell _s$.
\item[(c)]
$\alpha _{s,p}^{\vee}(\alpha _{t,q})=0$ and $\alpha _{s,p}^{\vee}(\alpha _{i_t})=0$ 
if $(t,q)<(s,p-1)$ in lexicographic order. 
\item[(d)]
$\alpha _{i_s}^{\vee}(\alpha _{s-1, 1})=-1$ and $\alpha _{i_s}^{\vee}(\alpha _{t,q})=0$
if $(t,q)<(s-1, 1)$ in lexicographic order. 
\item[(e)]
$\alpha _{0, 1}^{\vee}(\lambda )=1$ and $\alpha _{s,p}^{\vee}(\lambda )=0$ 
for $0 \le s \le k-1$ and $1\le p\le \ell _s$ such that $(s,p)\neq (0, 1)$. 
\item[(f)]
If $s\neq 0$ or $w_0 \neq 1$, then $\alpha _{i_s}^{\vee}(\lambda )=0$ for $0\le s \le k-1$.
\item[(g)]
If $w_0 =1$, then $\alpha _{i_1} ^{\vee}(\lambda )=1$.
\end{itemize}
\end{lem}

\subsection{Description of Demazure modules}
This subsection is devoted to the proof of Proposition \ref{4.3.1} below. 
Before doing this, we fix some notation. Let $\mathfrak{g}_{re} \subset \mathfrak{g}$ 
be the Levi subalgebra corresponding to $I^{re}$; 
note that $\mathfrak{g}_{re}$ is an ordinary Kac--Moody algebra. 
For a (fixed) reduced expression 
$v=r_{j_l} \cdots r_{j_2} r_{j_1} \in \mathcal{W}_{re}$, with $j_l ,\ldots ,j_2 ,j_1 \in I^{re}$, 
we set $F_v ^{\mathbf{m}}:=f_{j_l} ^{m_l} \cdots f_{j_2} ^{m_2} f_{j_1} ^{m_1}$ and 
$\widetilde{F}_v ^{\mathbf{m}}:=\Tilde{f}_{j_l} ^{m_l} \cdots \Tilde{f}_{j_2} ^{m_2} \Tilde{f}_{j_1} ^{m_1}$. 
Also, for $\mathbf{i}=(i_l ,\ldots , i_2 , i_1)$, set 
$F_{\mathbf{i}}^{\mathbf{m}}:=f_{i_l} ^{m_l} \cdots f_{i_2} ^{m_2} f_{i_1} ^{m_1}$ and 
$F_{\mathbf{i}}^{\max}:=f_{i_l} ^{\max}\cdots f_{i_1}^{\max}$, where 
$f_j ^{\max} u$, $j\in I^{re}$, $u\in V(\lambda)$, denotes the element 
$f_j ^p u \neq 0$, $p\ge 0$, such that $f_j ^{p+1} u$=0.
By convention, we write $F_v := F_v ^{\mathbf{m}}$ if $\mathbf{m}=(1,\ldots ,1 ,1)$. 
\begin{prop}\label{4.3.1}
Let $w\in \mathcal{W}$ satisfy conditions (\ref{eq4}) and (\ref{eq5}) for $\lambda \in P^+$, 
and fix a minimal dominant reduced expression $w=w_k r_{i_k}^{a_k} \cdots w_1 r_{i_1 }^{a_1} w_0$.
Then, with the notation above, we have
\begin{itemize}
\item[(1)]
$V_w (\lambda )
=\sum _{\underline{\mathbf{m}}, \underline{\bm{\epsilon }}} \mathbb{C}(q)F_{w_k}^{\mathbf{m}_k}
f_{i_k}^{\epsilon _k}\cdots F_{w_1}^{\mathbf{m}_1} f_{i_1}^{\epsilon _1}
F_{w_0}^{\mathbf{m}_0}u_{\lambda }$, where the summation is over all 
$\underline{\mathbf{m}}=(\mathbf{m}_s )_{s=0} ^k 
\in \prod _{s=0} ^k \mathbb{Z}_{\ge 0} ^{\ell (w_s)}$ 
and all $\underline{\bm{\epsilon }}=(\epsilon _t )_{t=1} ^k$, 
$0 \le \epsilon _t \le a_t$, $1\le t \le k$.
\item[(2)] $V_v (\lambda ) \subset V_w (\lambda)$ if $v \le w$ in $\mathcal{W}$.
\end{itemize}
\end{prop}
Part (2) of Proposition \ref{4.3.1} follows from part (1), and part (1) for $k=1$ 
is established by Lemmas \ref{4.3.2} and \ref{4.3.3} below; 
part (1) for $k\ge 2$ follows immediately by induction on $k$.

Now, suppose that $k=1$ for the expression of $w$ above, 
and write $w= R_{\mathbf{i}} r_i ^a R_{\mathbf{j}}$, 
with $a \in \mathbb{Z}_{\ge 1}$, $i\in I^{im},\ \mathbf{i}=(i_s ,\ldots ,i_2 ,i_1),\ 
\mathbf{j}=(j_t ,\ldots ,j_2 ,j_1)$, such that 
$R_{\mathbf{i}}, R_{\mathbf{j}}\in \mathcal{W}_{re}$ are reduced expressions.
In this case, we have $V_w (\lambda )=U_q ^+ (\mathfrak{g}) u_{w\lambda }$, 
where $u_{w\lambda }:=F_{\mathbf{i}}^{\max } f_i ^a F_{\mathbf{j}}u_{\lambda }
=F_{\mathbf{i}}^{\max } f_i ^a F_{\mathbf{j}} ^{\max} u_{\lambda }
\in V(\lambda )_{w\lambda }\setminus \{ 0\}$ by Lemma \ref{4.2.2}.
Set $J:=I\setminus \{ i\}$, and denote by $\mathfrak{g}_J \subset \mathfrak{g}$
the Levi subalgebra corresponding to $J$. If we write
$\bigl( e_i U_q ^+ (\mathfrak{g}_J) \bigl) ^m
:=\bigl( e_i U_q ^+ (\mathfrak{g}_J) \bigl)
\cdots \bigl( e_i U_q ^+ (\mathfrak{g}_J) \bigl)$ ($m$ times),
then $U_q ^+ (\mathfrak{g})=\bigoplus _{m=0}^{\infty}
U_q ^+ (\mathfrak{g}_J)\bigl( e_i U_q ^+ (\mathfrak{g}_J) \bigl) ^m$ 
and we have the following by weight consideration:
\begin{align}\label{eq7}
V_{w }(\lambda )=\bigoplus _{m=0} ^{\infty} U_q ^+ (\mathfrak{g}_J)
\bigl( e_i U_q ^+ (\mathfrak{g}_J ) \bigl) ^m u_{w \lambda }
=\bigoplus _{m=0} ^a U_q ^+ (\mathfrak{g}_J) \bigl( e_i U_q ^+ (\mathfrak{g}_J ) 
\bigl) ^m u_{w \lambda },
\end{align}
where the number $a$ in the right-hand side of (\ref{eq7}) is the one 
appearing in the expression $w =  R_{\mathbf{i}} r_i ^a R_{\mathbf{j}}$. 
By Lemma \ref{4.2.2}, we have
$ f_i ^m F_{\mathbf{j}} u_{\lambda }
\in V(\lambda ) _{r_i ^m R_{\mathbf{j}} \lambda } \setminus \{ 0 \}$ and 
$F_{\mathbf{i}} ^{\max }  f_i ^m F_{\mathbf{j}} u_{\lambda } \in 
V(\lambda )_{R_{\mathbf{i}} r_i ^m R_{\mathbf{j}} \lambda } \setminus \{ 0 \}$, $0\le m \le a$.
Note that $r_i ^m R_{\mathbf{j}} \lambda \in P^+$ for all $m\ge 1$ 
since $R_{\mathbf{i}}r_i ^m R_{\mathbf{j}}$ is a (minimal) dominant reduced expression.
Therefore, each $f_i ^m F_{\mathbf{j}} u_{\lambda }$ for 
$1\le m\le a$ is a $U_q (\mathfrak{g}_{re})$-highest weight vector.
\begin{lem}\label{4.3.2}
For each $1 \le m \le a,$ we have
$$U_q ^+ (\mathfrak{g}_J) \bigl( e_i U_q ^+ (\mathfrak{g}_J ) \bigl) ^{a-m}
u_{w \lambda } = \sum _{\mathbf{n} \in \mathbb{Z}_{\ge 0} ^s} \mathbb{C}(q)
F_{\mathbf{i}} ^{\mathbf{n}} f_i ^m F_{\mathbf{j}} u_{\lambda }.$$
\end{lem}
\textbf{Proof.}\ 
We proceed by descending induction on $m$.
If $m=a$, then we have 
$e_j U_q ^+ (\mathfrak{g}_J) u_{w\lambda} = \{ 0 \}$ 
for all $j \in J \setminus I^{re}$ by weight consideration.
Hence it follows that $U_q ^+ (\mathfrak{g}_J) u_{w \lambda }
=U_q ^+ (\mathfrak{g}_{re}) u_{w \lambda }$.
Also, $U_q ^+ (\mathfrak{g}_{re}) u_{w \lambda }$ 
is identical to 
$\sum _{\mathbf{n}} \mathbb{C}(q) F_{\mathbf{i}} ^{\mathbf{n}}  f_i ^a F_{\mathbf{j}} u_{\lambda }$
by [\textbf{Kas2}, Corollary 3.2.2]. Thus, Lemma \ref{4.3.2} follows in this case.

Suppose that $1\le m\le a-1$. By the induction hypothesis, we have
$$U_q ^+(\mathfrak{g}_J)
\bigl( e_i U_q ^+ (\mathfrak{g}_J ) \bigl) ^{a-m} u_{w\lambda }
=U_q ^+(\mathfrak{g}_J)e_i 
\sum _{\mathbf{n}} \mathbb{C}(q) F_{\mathbf{i}} ^{\mathbf{n}}f_i ^{m+1} F_{\mathbf{j}} u_{\lambda },$$
which is identical to 
$U_q ^+(\mathfrak{g}_J)
\sum \mathbb{C}(q) F_{\mathbf{i}} ^{\mathbf{n}} e_i f_i ^{m+1} F_{\mathbf{j}} u_{\lambda }$
since $i$ does not appear in $\mathbf{i}$ and each $F_{\mathbf{i}} ^{\mathbf{n}}$ 
in the summation commutes with $e_i $. Since $a_{ii}=0$ by Lemma \ref{4.2.1}, it follows that 
$$e_i f_i ^{m+1} F_{\mathbf{j}} u_{\lambda }
=\underbrace{\{ m+1 \} _i}_{=m+1(\ge 2)}
\Bigl[ \underbrace{\alpha _i ^{\vee}(\lambda -\alpha _{j_1}-
\cdots -\alpha _{j_t})}_{=1}-\underbrace{\frac{1}{2} a_{ii} m}_{=0} \Bigl] _i 
f_i ^m F_{\mathbf{j}} u_{\lambda }$$
by [\textbf{JKK}, Lemma 2.5]. Therefore, we obtain
\begin{align}\label{eq8}
U_q ^+(\mathfrak{g}_J) \sum \mathbb{C}(q) F_{\mathbf{i}} ^{\mathbf{n}} 
e_i f_i ^{m+1} F_{\mathbf{j}} u_{\lambda }
=U_q ^+(\mathfrak{g}_J)\sum \mathbb{C}(q) F_{\mathbf{i}} ^{\mathbf{n}} 
f_i ^m F_{\mathbf{j}} u_{\lambda }.
\end{align}
Also, we see that the right-hand side of (\ref{eq8}) is identical to 
$U_q ^+(\mathfrak{g}_{re})\sum \mathbb{C}(q) F_{\mathbf{i}} ^{\mathbf{n}} 
f_i ^m F_{\mathbf{j}} u_{\lambda }$ by weight consideration, 
which is identical to
$\sum \mathbb{C}(q) F_{\mathbf{i}} ^{\mathbf{n}} f_i ^m F_{\mathbf{j}} u_{\lambda }$ 
by [\textbf{Kas2}, Corollary 3.2.2]. This proves Lemma \ref{4.3.2}. \qed
\begin{lem}\label{4.3.3}
We have
$U_q ^+ (\mathfrak{g}_J) \bigl( e_i U_q ^+ (\mathfrak{g}_J ) \bigl) ^a
u_{w \lambda } = \sum _{\mathbf{n}, \mathbf{m}} \mathbb{C}(q)
F_{\mathbf{i}} ^{\mathbf{n}} F_{\mathbf{j}} ^{\mathbf{m}} u_{\lambda }.$
\end{lem}
\textbf{Proof.}\ 
If we set $m=0$ in the proof of Lemma \ref{4.3.2}, then the same argument shows that
\begin{align}\label{eq9}
U_q ^+ (\mathfrak{g}_J)\bigl( e_i U_q ^+ (\mathfrak{g}_J)\bigl) ^a
u_{w\lambda } =U_q ^+ (\mathfrak{g}_{re})
\sum \mathbb{C}(q) F_{\mathbf{i}} ^{\mathbf{n}} F_{\mathbf{j}} ^{\mathbf{m}} u_{\lambda }.
\end{align}
Now, it is easily seen that the right-hand side of (\ref{eq9}) is identical to 
$\sum \mathbb{C}(q) F_{\mathbf{i}} ^{\mathbf{n}} F_{\mathbf{j}} ^{\mathbf{m}} u_{\lambda }$ 
by direct calculation. This proves Lemma \ref{4.3.3}. \qed

\subsection{The maps $\Psi _{\lambda ,w}$ and $\Phi _{\lambda ,w}$}
Let $\Tilde{\mathfrak{g}}_{re} \subset \Tilde{\mathfrak{g}}$ denote the Levi subalgebra 
corresponding to $\Tilde{I}^{re}:=\{ (i,1) \} _{i\in I^{re}}$. 
Since $\mathfrak{g}_{re}$ (see \S 4.3) and $\Tilde{\mathfrak{g}}_{re}$
are ordinary Kac--Moody algebras with the same Cartan datum, 
there exists an embedding of $\mathbb{C}(q)$-algebras:
\begin{align}\label{eq10}
\Psi _{re} :U_q (\mathfrak{g}_{re})\hookrightarrow 
U_q (\Tilde{\mathfrak{g}}_{re}),\ 
e_i \longmapsto e_{(i,1)},\ 
f_i \longmapsto f_{(i,1)},\ 
q^{\alpha _i ^{\vee}} \longmapsto q^{\Tilde{\alpha}_{(i,1)} ^{\vee}}.
\end{align}
According to the minimal dominant reduced expression
$w=w_k r_{i_k}^{a_k} \cdots w_1 r_{i_1 }^{a_1} w_0$, 
the element $\Tilde{w} \in \widetilde{W}$ (see Remark \ref{4.1.2}) 
is expressed as $\Tilde{w} = \Tilde{w}_k R_{(i_k, \bm{a}_k)} \cdots 
\Tilde{w}_1 R_{(i_1, \bm{a}_1)} \Tilde{w}_0$, with $\bm{a}_s =(a_s , \ldots ,2,1),\ 1\le s \le k$, 
where each $\Tilde{w}_s \in \widetilde{W}_{re}$ corresponds to 
$w_s \in \mathcal{W}_{re}$ via the isomorphism (\ref{eq3}) in \S 3.2 
and we write $(i,\mathbf{m})=((i, m_l), \ldots , (i, m_1))$ if $\mathbf{m}=(m_l ,\ldots , m_1)$; 
recall the notation $R_{\mathbf{(i,m)}} =r_{(i_l , m_l)} \cdots r_{(i_1 , m_1)}$.
Note that this expression for $\Tilde{w}$ has properties similar 
to those in Lemma \ref{4.2.2} for $\Tilde{\lambda} \in \widetilde{P}^+$.
Let $\widetilde{V}_{\Tilde{w}} (\Tilde{\lambda})$ denote the 
Demazure submodule, corresponding to $\Tilde{w}$, of the irreducible highest weight 
$U_q (\Tilde{\mathfrak{g}})$-module $\widetilde{V}(\Tilde{\lambda})$
of highest weight $\Tilde{\lambda}$. We can think of $\widetilde{V}_{\Tilde{w}} (\Tilde{\lambda})$
as a $U_q ^+ (\mathfrak{g}_{re})$-module via $\Psi _{re}$ (see (\ref{eq10})). 

In this subsection, we construct a $U_q ^+ (\mathfrak{g}_{re})$-linear maps 
$\Psi _{\lambda ,w}: V_w (\lambda ) \hookrightarrow \widetilde{V}_{\Tilde{w}}(\Tilde{\lambda })$ and 
$\Phi _{\lambda ,w}: \widetilde{V}_{\Tilde{w}}(\Tilde{\lambda})\twoheadrightarrow V_w (\lambda)$ 
such that $\Phi _{\lambda ,w} \circ \Psi _{\lambda ,w} = \mathrm{id}_{V_w (\lambda)}$. 
For this purpose, we first study the $U_q ^+ (\mathfrak{g}_{re})$-module structure of 
$V_w (\lambda )$ and $\widetilde{V}_{\Tilde{w}}(\Tilde{\lambda})$.
As in \S 4.3, define $F_{\Tilde{v}} ^{\mathbf{m}}, \widetilde{F}_{\Tilde{v}} ^{\mathbf{m}}$ and 
$F_{(\mathbf{i}, \mathbf{n})} ^{\mathbf{m}}$ for $\Tilde{v} \in  \widetilde{W}$ and 
$(\mathbf{i,n}) \in \widetilde{\mathcal{I}}$. 
Set $F^{\underline{\mathbf{m}}} _{\Tilde{w}_k ,\ldots ,\Tilde{w}_s}:=
F^{\mathbf{m}_k}_{\Tilde{w}_k} \cdots F^{\mathbf{m}_s}_{\Tilde{w}_s}
\in U_q ^- (\Tilde{\mathfrak{g}})$, $\underline{\mathbf{m}} = (\mathbf{m}_t) _{t=s} ^k
\in \prod _{t=s} ^k \mathbb{Z}_{\ge 0} ^{\ell (\Tilde{w}_t)}$, and 
$\Tilde{u}_{\underline{\mathbf{n}}}:=F_{(i_s , \mathbf{n}_s)} F_{\Tilde{w}_{s-1}}
\cdots F_{(i_1 , \mathbf{n}_1)} F_{\Tilde{w}_0}\Tilde{u}_{\Tilde{\lambda}}\in \widetilde{V}(\Tilde{\lambda})$, 
$\underline{\mathbf{n}} = (\mathbf{n}_t) _{t=1} ^s$, 
$\mathbf{n}_t \subset \{ 1,2, \ldots , a_t \}$. 
Here we understand that $\Tilde{u}_{\underline{\mathbf{n}}}=\Tilde{u}_{\Tilde{\lambda}}$ if $s=0$.
Note that if $\# \mathbf{n}_t >1$, then $a_{i_t , i_t}=0$ by Lemma \ref{4.2.1}, 
and hence all the $f_{(i_t , m)}$, $m\ge 1$, commute with each other. 
Therefore, the element $F_{(i_t , \mathbf{n}_t)}$ is independent of 
the choice of an ordering of the elements in $\mathbf{n}_t $.
Also, set $F^{\underline{\mathbf{m}}} _{w_k ,\ldots , w_s}:=
F^{\mathbf{m}_k}_{w_k} \cdots F^{\mathbf{m}_s}_{w_s} \in U_q ^- (\mathfrak{g})$ and 
$u_{\underline{\mathbf{n}}}:=f_{i_s} ^{\# \mathbf{n}_s} F_{w_{s-1}}
\cdots f_{i_1} ^{\# \mathbf{n}_1} F_{w_0} u_{\lambda} \in V(\lambda )$; 
note that $u_{\underline{\mathbf{n}}}$ depends only on the cardinalities 
$\# \mathbf{n}_t$, $1\le t \le s$. Hence we write
$u_{\underline{\bm{\epsilon }}} := u_{\underline{\mathbf{n}}}$
for $\underline{\bm{\epsilon }}=(\epsilon _t) _{t=1} ^s$
if $\epsilon _t = \# \mathbf{n}_t$, $1 \le t \le s$.
If $\mathbf{n}_t \neq \emptyset$ for $1\le t \le s$, then 
$u_{\underline{\mathbf{n}}}$ (resp., $\Tilde{u}_{\underline{\mathbf{n}}}$)
is a highest weight vector for $U_q (\mathfrak{g}_{re})$ (resp., $U_q (\Tilde{\mathfrak{g}}_{re})$)
of weight $r_{i_s} ^{\# \mathbf{n}_s} w_{s-1} \cdots r_{i_1} ^{\# \mathbf{n}_1} w_0 (\lambda) \in P$ 
(resp., $R_{(i_s ,\mathbf{n}_s)} \Tilde{w}_{s-1} \cdots R_{(i_1 , \mathbf{n}_1)} \Tilde{w}_0 (\Tilde{\lambda})
\in \widetilde{P}$) by Lemma \ref{4.2.2}. Moreover, in this case, we have 
$\alpha _i ^{\vee } \bigl( r_{i_s} ^{\# \mathbf{n}_s} w_{s-1} \cdots r_{i_1} ^{\# \mathbf{n}_1} w_0 (\lambda) \bigl) =
\Tilde{\alpha }_{(i,1)} ^{\vee} \bigl( R_{(i_s ,\mathbf{n}_s)} \Tilde{w}_{s-1} \cdots R_{(i_1 , \mathbf{n}_1)} \Tilde{w}_0 (\Tilde{\lambda}) \bigl)$ for all $i\in I^{re}$, and hence an isomorphism 
$U_q (\mathfrak{g}_{re}) u_{\underline{\mathbf{n}}} \xrightarrow{\cong} 
U_q (\Tilde{\mathfrak{g}}_{re}) \Tilde{u}_{\underline{\mathbf{n}}}$
of $U_q (\mathfrak{g}_{re})$-modules via $\Psi _{re}$ (see (\ref{eq10})).
\begin{lem}\label{4.4.1}
As $U_q ^+(\mathfrak{g}_{re})$-modules, $V_w (\lambda)$ and 
$\widetilde{V}_{\Tilde{w}}(\Tilde{\lambda})$ decompose as follows: 
\begin{itemize}
\item[(1)]
$V_{w}(\lambda )=
\bigoplus _{s=0} ^k \bigoplus _{\underline{\bm{\epsilon}}}\bigl( \sum _{\underline{\mathbf{m}}}
\mathbb{C}(q) F^{\underline{\mathbf{m}}} _{w_k ,\ldots , w_s} u_{\underline{\bm{\epsilon}}}\bigl)$, 
where the summation is over all 
$\underline{\mathbf{m}} = (\mathbf{m}_t ) _{t=s} ^k
\in \prod _{t=s} ^k \mathbb{Z}_{\ge 0} ^{\ell (w_t)}$ and all 
$\underline{\bm{\epsilon}}=(\epsilon _t )_{t=1} ^s ,$ with 
$1\le \epsilon _t \le a_t ,$ $1\le t \le s ,$ for each $0\le s \le k$.
\item[(2)]
$\widetilde{V}_{\Tilde{w}}(\Tilde{\lambda}) =
\bigoplus _{s=0} ^k \bigoplus _{\underline{\mathbf{n}}}\bigl( \sum _{\underline{\mathbf{m}}}
\mathbb{C}(q) F^{\underline{\mathbf{m}}} _{\Tilde{w}_k ,\ldots ,\Tilde{w}_s}
\Tilde{u}_{\underline{\mathbf{n}}} \bigl)$, 
where the summation is over all 
$\underline{\mathbf{m}} = (\mathbf{m}_t) _{t=s} ^k
\in \prod _{t=s} ^k \mathbb{Z}_{\ge 0} ^{\ell (\Tilde{w}_t)}$ and all 
$\underline{\mathbf{n}}=(\mathbf{n}_t )_{t=1} ^s ,$ with 
$\emptyset \neq \mathbf{n}_t \subset \{ 1,2, \ldots , a_t \} ,$ $1\le t \le s ,$ 
for each $0\le s \le k$.
\end{itemize}
\end{lem}
\textbf{Proof.}
We give a proof only for (1); the proof for (2) is similar.
By Lemma \ref{4.2.2} and Proposition \ref{4.3.1} (1), 
it is easily seen that $\{ u_{\underline{\bm{\epsilon }}} \} _{\underline{\bm{\epsilon }}}$, 
where the $\underline{\bm{\epsilon }}$ runs as in Lemma \ref{4.4.1} (1), forms a 
complete set of representatives of linearly independent 
$U_q (\mathfrak{g}_{re})$-highest weight vectors in $V _w (\lambda )$. 
Therefore, the summation 
$\sum _{s=0} ^k \sum _{\underline{\bm{\epsilon }}} \bigl( \sum _{\underline{\mathbf{m}}}
\mathbb{C}(q) F^{\underline{\mathbf{m}}} _{w _k ,\ldots , w_s } u_{\underline{\bm{\epsilon }}} \bigl)$
is a direct sum with respect to $s$ and $\underline{\bm{\epsilon }}$. 
Now, the equality in (1) follows from Proposition \ref{4.3.1} (1). 
Also, from the arguments in \S 4.3, we know that each direct summand 
$\sum _{\underline{\mathbf{m}}}
\mathbb{C}(q) F^{\underline{\mathbf{m}}} _{w_k ,\ldots , w_s} u_{\underline{\bm{\epsilon}}}$ 
is stable under the action of $U_q ^+(\mathfrak{g}_{re})$.
Thus, the expression for $V_w (\lambda )$ in Lemma \ref{4.4.1} (1) 
is a decomposition as $U_q ^+ (\mathfrak{g}_{re})$-modules. \qed

\vspace{2mm}
By using the descriptions of $V_w (\lambda )$ and 
$\widetilde{V}_{\Tilde{w}}(\Tilde{\lambda})$ in Lemma \ref{4.4.1}, define 
$\Phi _{\lambda ,w}: \widetilde{V}_{\Tilde{w}}(\Tilde{\lambda})\twoheadrightarrow V_w (\lambda)$ 
by $F^{\underline{\mathbf{m}}} _{\Tilde{w}_k ,\ldots ,\Tilde{w}_s}
\Tilde{u}_{\underline{\mathbf{n}}} \mapsto 
F^{\underline{\mathbf{m}}} _{w_k ,\ldots , w_s}
u_{\underline{\mathbf{n}}}$, and 
$\Psi _{\lambda , w}: V_w (\lambda ) \hookrightarrow \widetilde{V}_{\Tilde{w}}(\Tilde{\lambda})$ by
$F^{\underline{\mathbf{m}}} _{w_k ,\ldots , w_s} u_{\underline{\bm{\epsilon}}}
=F^{\underline{\mathbf{m}}} _{w_k ,\ldots , w_s}
u_{\underline{\mathbf{n}}} \mapsto 
F^{\underline{\mathbf{m}}} _{\Tilde{w}_k ,\ldots ,\Tilde{w}_s}
\Tilde{u}_{\underline{\mathbf{n}}}$ for 
$\underline{\mathbf{n}}=(\mathbf{n}_t)_{t= 1} ^s$, 
$\mathbf{n}_t =\{ 1,2,\ldots , \epsilon _t \}$, 
$1 \le \epsilon _t \le a_t$, $1\le t \le s$, with 
$\underline{\bm{\epsilon }}=(\epsilon _t )_{t=1} ^s$.
Since $U_q (\Tilde{\mathfrak{g}}_{re}) \Tilde{u}_{\underline{\mathbf{n}}}
\cong U_q (\mathfrak{g}_{re}) u_{\underline{\mathbf{n}}}$ as $U_q (\mathfrak{g}_{re})$-modules, 
we deduce that $\sum _{\underline{\mathbf{m}}}
\mathbb{C}(q) F^{\underline{\mathbf{m}}} _{\Tilde{w}_k ,\ldots ,\Tilde{w}_s}
\Tilde{u}_{\underline{\mathbf{n}}} \xrightarrow{\cong}
\sum _{\underline{\mathbf{m}}}
\mathbb{C}(q) F^{\underline{\mathbf{m}}} _{w_k ,\ldots ,w_s}
u_{\underline{\mathbf{n}}}$ as $U_q ^+ (\mathfrak{g}_{re})$-modules
via the map $\Phi _{\lambda ,w}$ for all $\underline{\mathbf{n}}$ as in 
Lemma \ref{4.4.1} (2). Thus, $\Phi _{\lambda ,w}$ is 
well-defined, surjective, and $U_q ^+ (\mathfrak{g}_{re})$-linear.
Also, we can verify that $\Psi _{\lambda ,w}$ is well-defined, injective, 
and $U_q ^+ (\mathfrak{g}_{re})$-linear. Note that $\Psi _{\lambda ,w}$ is also surjective 
(and hence an isomorphism) if and only if $a_1 =\cdots =a_k =1$.
By the definitions, we have $\Phi _{\lambda ,w} \circ \Psi _{\lambda ,w} = \mathrm{id}_{V_w (\lambda)}$.

\subsection{Global bases and diagram automorphisms}
In this subsection, we recall some fundamental properties of 
the global basis of $\widetilde{V} _{\Tilde{w}} (\Tilde{\lambda })$. 
Also, we study certain symmetries of $\widetilde{V} _{\Tilde{w}} (\Tilde{\lambda })$, 
which come from diagram automorphisms of $\Tilde{\mathfrak{g}}$.

Let $\widetilde{B}(\Tilde{\lambda })$ and 
$\{ \widetilde{G}_{\Tilde{\lambda }}(b')\} _{ b' \in \widetilde{B}(\Tilde{\lambda })}$ 
denote the crystal basis and the global basis of 
$\widetilde{V}(\Tilde{\lambda })$ with the crystal lattice 
$\widetilde{L}(\Tilde{\lambda})$, and let 
$\widetilde{V}(\Tilde{\lambda })^{\textbf{A}} \subset \widetilde{V}(\Tilde{\lambda })$ 
denote the $\mathbf{A}$-form of $\widetilde{V}(\Tilde{\lambda })$. 
If we set 
\begin{align*}
\widetilde{B}_{\Tilde{w}}(\Tilde{\lambda }):=\{ & 
\widetilde{F}_{\Tilde{w}_k}^{\mathbf{m}_k}
\widetilde{F}_{(i_k , \mathbf{n}_k)} \cdots 
\widetilde{F}_{\Tilde{w}_1}^{\mathbf{m}_1}
\widetilde{F}_{(i_1 , \mathbf{n}_1)}
\widetilde{F}_{\Tilde{w}_0}^{\mathbf{m}_0}
\Tilde{u}_{\Tilde{\lambda}}
\mod q\widetilde{L}(\Tilde{\lambda}) \mid \\
& (\mathbf{m}_s)_{s=0}^k \in \prod _{s=0} ^k \mathbb{Z}_{\ge 0} ^{\ell (\Tilde{w}_s)},
\ \mathbf{n}_t \subset \{ 1,2,\ldots ,a_t \},\ 1\le t \le k \} \setminus \{ 0 \}
\subset \widetilde{B}(\Tilde{\lambda}), 
\end{align*}
then we have
$\widetilde{V}_{\Tilde{w}}(\Tilde{\lambda })
=\bigoplus _{b' \in \widetilde{B}_{\Tilde{w}}(\Tilde{\lambda })}
\mathbb{C}(q)\widetilde{G}_{\Tilde{\lambda }}(b')$
by [\textbf{Kas2}, Proposition 3.2.3].
Let $B_w (\lambda ) \subset B(\lambda )$ denote the inverse image of 
$\widetilde{B}_{\Tilde{w}}(\Tilde{\lambda })$ under the embedding 
$B(\lambda )\hookrightarrow \widetilde{B}(\Tilde{\lambda })$ 
of (\ref{eq2}) in \S 2.5, and 
$\widetilde{B}_{\Tilde{w}}^0 (\Tilde{\lambda }) \subset \widetilde{B}(\Tilde{\lambda})$ 
the image of $B_w (\lambda )$ under this embedding. Then, these subsets can be written as 
\begin{align*}
B_w (\lambda ) = \{ & 
\widetilde{F}_{w_k} ^{\mathbf{m}_k} 
\Tilde{f}_{i_k}^{\epsilon _k} \cdots 
\widetilde{F}_{w_1} ^{\mathbf{m}_1} 
\Tilde{f}_{i_1} ^{\epsilon _1} 
\widetilde{F}_{w_0} ^{\mathbf{m}_0} u_{\lambda} 
\mod qL(\lambda) \mid \\
& (\mathbf{m}_s )_{s=0} ^k \in 
\prod _{s=0} ^k \mathbb{Z}_{\ge 0} ^{\ell (w_s)},\ 
0 \le \epsilon _t \le a_t ,\ 1\le t \le k \} \setminus \{ 0 \} ,\ \text{and} \\ 
\widetilde{B}_{\Tilde{w}} ^0 (\Tilde{\lambda }) = \{ & 
\widetilde{F}_{\Tilde{w}_k}^{\mathbf{m}_k}
\widetilde{F}_{(i_k , \mathbf{n}_k)} \cdots 
\widetilde{F}_{\Tilde{w}_1}^{\mathbf{m}_1}
\widetilde{F}_{(i_1 , \mathbf{n}_1)}
\widetilde{F}_{\Tilde{w}_0}^{\mathbf{m}_0}
\Tilde{u}_{\Tilde{\lambda}}
\mod q\widetilde{L}(\Tilde{\lambda}) \mid \\
& (\mathbf{m}_s)_{s=0}^k \in \prod _{s=0} ^k \mathbb{Z}_{\ge 0} ^{\ell (\Tilde{w}_s)},
\ \mathbf{n}_t = \{ 1,2,\ldots , \epsilon _t \} ,\ 0 \le \epsilon _t \le a_t , \ 
1\le t \le k \} \setminus \{ 0 \} .
\end{align*}
In \S 5, we will show that the subset $B_w (\lambda ) \subset B(\lambda )$ 
satisfies the condition in Theorem \ref{thm3}.

Let us set $\Tilde{J}_s :=\{ (i_s ,m) \in \Tilde{I} \mid 1\le m\le a_s \}$ for $1 \le s \le k$, 
with each $i_s \in I^{im}$ and $a_s$ appearing in the expression 
$w=w_k r_{i_k} ^{a_k } \cdots w_1 r_{i_1} ^{a_1} w_0$.
Also, we set 
$\Omega _{\lambda ,w} :=\prod _{s=1}^k \mathfrak{S}(\Tilde{J}_s) \subset \Omega $, 
where $\mathfrak{S}(\Tilde{J}_s)$ denote the permutation group on $\Tilde{J}_s$ (see \S 2.2). 
Note that the action of $\omega \in \Omega _{\lambda ,w}$ on $\Tilde{I}$ is expressed as: 
$\omega (j,n) =(j,n)$ if $(j,n) \notin \bigcup _{s=1}^k \Tilde{J}_s$, and 
$\omega (i_s , m) = (i_s , m')$ for some $1\le m' \le a_s$ if $(i_s , m) \in \Tilde{J}_s$.
From this observation, we can show that $\widetilde{V}_{\Tilde{w}} (\Tilde{\lambda })$ 
and $\widetilde{B}_{\Tilde{w}} (\Tilde{\lambda })$ is stable under the action of $\Omega _{\lambda ,w}$.
By using the description of Lemma \ref{4.4.1} for $\widetilde{V}_{\Tilde{w}} (\Tilde{\lambda })$, 
each $\omega \in \Omega _{\lambda ,w}$ sends the element 
$F^{\underline{\mathbf{m}}} _{\Tilde{w}_k ,\ldots ,\Tilde{w}_s}
\Tilde{u}_{\underline{\mathbf{n}}}$, with $\underline{\mathbf{n}}=(\mathbf{n}_s ,\ldots , \mathbf{n}_1)$, 
to the element of the form 
$F^{\underline{\mathbf{m}}} _{\Tilde{w}_k ,\ldots ,\Tilde{w}_s}
\Tilde{u}_{\underline{\mathbf{n}'}}$, with $\underline{\mathbf{n}'}=(\mathbf{n}' _s ,\ldots , \mathbf{n}' _1)$ 
such that $\# \mathbf{n}_t = \# \mathbf{n}' _t$ for $1\le t \le s$.
Note that this expression for the action of $\Omega _{\lambda ,w}$ on 
$\widetilde{V}_{\Tilde{w}} (\Tilde{\lambda })$ shows that 
$\Phi _{\lambda ,w} \circ \omega = \Phi _{\lambda ,w}$ for all $\omega \in \Omega _{\lambda ,w}$.
Also, for each $\underline{\mathbf{n}'}=(\mathbf{n}' _s ,\ldots , \mathbf{n}' _1)$
with $\# \mathbf{n}' _t = \# \mathbf{n}_t$ for $1\le t \le s$, 
we can find $\omega \in \Omega _{\lambda ,w}$ such that 
$\omega (F^{\underline{\mathbf{m}}} _{\Tilde{w}_k ,\ldots ,\Tilde{w}_s}
\Tilde{u}_{\underline{\mathbf{n}}}) = 
F^{\underline{\mathbf{m}}} _{\Tilde{w}_k ,\ldots ,\Tilde{w}_s}\Tilde{u}_{\underline{\mathbf{n}'}}$.
In particular, if we take $\underline{\mathbf{n}'}=(\mathbf{n}' _s ,\ldots , \mathbf{n}' _1)$
such that $\mathbf{n}' _t = \{1,2, \ldots , \# \mathbf{n} _t \}$ for $1\le t \le s$, then 
the corresponding $\omega \in \Omega _{\lambda ,w}$ (not necessarily unique) sends 
$F^{\underline{\mathbf{m}}} _{\Tilde{w}_k ,\ldots ,\Tilde{w}_s}\Tilde{u}_{\underline{\mathbf{n}}}$
into the image of the map $\Psi _{\lambda ,w}$. 
By the same argument, for each $b' \in \widetilde{B}_{\Tilde{w}} (\Tilde{\lambda })$, 
we can find $\omega \in \Omega _{\lambda ,w}$  
such that $\omega (b')$ lies in $\widetilde{B}_{\Tilde{w}} ^0 (\Tilde{\lambda })$. Namely, we have 
$\widetilde{B}_{\Tilde{w}} (\Tilde{\lambda })=
\Omega _{\lambda ,w} \widetilde{B}_{\Tilde{w}} ^0 (\Tilde{\lambda })$.
Note that 
$\widetilde{G}_{\Tilde{\lambda }} \circ \omega 
=\omega \circ \widetilde{G}_{\Tilde{\lambda }}$ for all $\omega \in \Omega _{\lambda ,w}$
(see [\textbf{S}, Lemma 3.4]).

\subsection{Compatibility with the global bases}
This subsection is devoted to the proof of Proposition \ref{4.6.1} below.
For $b' \in \widetilde{B}_{\Tilde{w}}(\Tilde{\lambda})$, we set 
$\mathbb{G}(b'):= \Phi _{\lambda ,w}\bigl( \widetilde{G}_{\Tilde{\lambda}}(b')\bigl) \in V_w (\lambda )$.
\begin{prop}\label{4.6.1}
The element $\mathbb{G}(b') \in V_w (\lambda),$ 
$b' \in \widetilde{B}_{\Tilde{w}}(\Tilde{\lambda})$, is a global basis element.
\end{prop}
\textbf{Proof.}
To complete the proof, in view of the characterization 
of a global basis element (see \S 2.1), we must show 
(i) $\overline{\mathbb{G}(b')} = \mathbb{G}(b')$, 
(ii) $\mathbb{G}(b') \in V(\lambda)^{\mathbf{A}} \cap L(\lambda)$, and 
(iii) $\mathbb{G}(b') \mod qL(\lambda) \in B(\lambda)\setminus \{ 0\}$. 
Note that we know from \S 4.5 that 
(i)$'$ $\overline{\widetilde{G}_{\Tilde{\lambda }}(b')}=\widetilde{G}_{\Tilde{\lambda }}(b')$, 
(ii)$'$ $\widetilde{G}_{\Tilde{\lambda }}(b')\in \widetilde{V}_{\Tilde{w}}(\Tilde{\lambda })
\cap \widetilde{V}(\Tilde{\lambda })^{\textbf{A}} \cap \widetilde{L}(\Tilde{\lambda })$, and 
(iii)$'$ $\widetilde{G}_{\Tilde{\lambda }}(b')\equiv b' \mod q\widetilde{L}(\Tilde{\lambda })$. 
Then, the equality (i) follows from (i)$'$. 
Obviously, we have 
$\Phi _{\lambda ,w}\bigl( \widetilde{V}_{\Tilde{w}}(\Tilde{\lambda })
\cap \widetilde{V}(\Tilde{\lambda })^{\textbf{A}}\bigl) \subset V(\lambda )^{\textbf{A}}$.
Therefore, by (ii)$'$, it suffices to verify that 
$\Phi _{\lambda ,w}\bigl( \widetilde{V}_{\Tilde{w}}(\Tilde{\lambda })
\cap \widetilde{L}(\Tilde{\lambda})\bigl) \subset L(\lambda)$ to show (ii).
For this purpose, the following is enough:
\begin{claim}\label{cl1}
If $\widetilde{F}_{(\mathbf{i}, \mathbf{m})}$ is of the form 
$\widetilde{F}_{\Tilde{w}_k}^{\mathbf{m}_k}
\widetilde{F}_{(i_k , \mathbf{n}_k)} \cdots 
\widetilde{F}_{\Tilde{w}_1}^{\mathbf{m}_1}
\widetilde{F}_{(i_k , \mathbf{n}_1)}
\widetilde{F}_{\Tilde{w}_0}^{\mathbf{m}_0}$, with
$(\mathbf{m}_s)_{s=0}^k \in \prod _{s=0} ^k \mathbb{Z}_{\ge 0} ^{\ell (\Tilde{w}_s)},
\ \mathbf{n}_t \subset \{ 1,2,\ldots ,a_t \},\ 1\le t \le k$, 
then we have $\Phi _{\lambda ,w}(\widetilde{F}_{(\mathbf{i}, \mathbf{m})} \Tilde{u}_{\Tilde{\lambda}})=\widetilde{F}_{\mathbf{i}}u_{\lambda}$. 
\end{claim}
\textbf{Proof of Claim 1.} 
Set $v = \widetilde{F}_{(\mathbf{i}, \mathbf{m})} \Tilde{u}_{\Tilde{\lambda}}$. 
We proceed by induction on the length $l$ of $\mathbf{(i,m)}=((i_s, m_s))_{s=1} ^l$.
We assume that $\Phi _{\lambda ,w}(v)=\widetilde{F}_{\mathbf{i}}u_{\lambda}$, and show that
$\Phi _{\lambda ,w}(\Tilde{f}_{(i,m)} v)=\Tilde{f}_i \widetilde{F}_{\mathbf{i}}u_{\lambda}$
for all $(i,m) \in \Tilde{I}$ such that $\Tilde{f}_{(i,m)} \widetilde{F}_{(\mathbf{i,m})}$
is also of the form given as in Claim 1.
\begin{flushleft}
\underline{Case 1:\ $i\in I^{im}$.}
\end{flushleft}
Since $(i,m)$ does not appear in $\mathbf{(i,m)}$, we have
$v \in \mathrm{Ker}(e_{(i,m)})$ and hence $\Tilde{f}_{(i,m)}v =f_{(i,m)}v$. Hence we deduce that 
$\Phi _{\lambda ,w}(\Tilde{f}_{(i,m)}v)= \Phi _{\lambda ,w}(f_{(i,m)}v) =f_i \Phi _{\lambda ,w}(v)$.
Also, we have $\Tilde{f}_i \Phi _{\lambda ,w}(v)=f_i \Phi _{\lambda ,w}(v)$
since $\Tilde{f}_i = f_i$ if $i\in I^{im}$ (see \S 2.1).
Thus the assertion follows in this case.
\begin{flushleft}
\underline{Case 2: $i\in I^{re}$.}
\end{flushleft}
Let $v=\sum _{n \ge 0} f_{(i,1)} ^{(n)} v_n$, with 
$v_n \in \mathrm{Ker}(e_{(i,1)})$, be the $(i,1)$-string decomposition (see \S 2.1).
Since $\widetilde{V}_{\Tilde{w}}(\Tilde{\lambda })$ is stable under the action of $e_{(i,1)}$, 
we can show, by induction on $n$, that 
$v_n$, $f_{(i,1)} ^{(n)} v_n \in \widetilde{V}_{\Tilde{w}}(\Tilde{\lambda })$ for all $n\ge 0$.
Hence the equality $\Phi _{\lambda ,w}(v)=\sum _{n \ge 0} f_i ^{(n)} \Phi _{\lambda ,w}(v_n)$ holds.
Also, we have $\Phi _{\lambda ,w}(v_n) \in \mathrm{Ker}(e_i)$ since 
$e_i \Phi _{\lambda ,w}(v_n )=\Phi _{\lambda ,w }(e_{(i,1)} v_n) =0$. 
Therefore, this expression for $\Phi _{\lambda ,w}(v)$ gives the $i$-string decomposition.
Thus we have $\Tilde{f} _i \Phi _{\lambda ,w}(v)=\sum _{n \ge 0} f_i ^{(n+1)} \Phi _{\lambda ,w}(v_n)$, 
which is also identical to the image of the element 
$\Tilde{f}_{(i,1)} v = \sum _{n\ge 0} f_{(i,1)} ^{(n+1)} v_n$
under the map $\Phi _{\lambda ,w}$. This proves Claim \ref{cl1}. \bqed

\vspace{2mm}
Now, let us show (iii). 
Note that Claim \ref{cl1} above shows that $\Phi _{\lambda ,w}$ induces the map 
$\overline{\Phi} _{\lambda ,w}: \widetilde{B}_{\Tilde{w}} (\Tilde{\lambda}) \rightarrow B(\lambda)$. 
Moreover, it is easily seen that the image of this map is 
$B_w (\lambda )$, and that the equality 
$\overline{\Phi} _{\lambda ,w} \circ \widetilde{\ \ }=\mathrm{id}_{B_w (\lambda)}$ holds, 
where $\widetilde{\ \ }$ denotes (the restriction of) the map of (\ref{eq2}) in \S 2.5.
From these observations, we see that if we take an element 
$\Tilde{b} \in \widetilde{B}_{\Tilde{w}}^0 (\Tilde{\lambda })$,
which is the image of an element $b\in B_w(\lambda )$ under the map $\widetilde{\ \ }$, then we have 
$\mathbb{G}(\Tilde{b})=\Phi _{\lambda ,w}\bigl( \widetilde{G}_{\Tilde{\lambda }}(\Tilde{b})\bigl)
\equiv  \overline{\Phi} _{\lambda ,w}(\Tilde{b}) \equiv b \mod qL(\lambda )$, 
where the second equality is due to (iii)$'$. 
For a general element $b'\in \widetilde{B}_{\Tilde{w}}(\Tilde{\lambda })$,
there exists a diagram automorphism $\omega \in \Omega _{\lambda ,w}$ such that 
$\omega ^{-1}(b')=\Tilde{b} \in \widetilde{B}_{\Tilde{w}} ^0(\Tilde{\lambda })$, where $\Tilde{b}$ is as above. 
Since $\widetilde{G}_{\Tilde{\lambda }} \circ \omega = \omega \circ \widetilde{G}_{\Tilde{\lambda }}$ and 
$\Phi _{\lambda ,w} \circ \omega = \Phi _{\lambda ,w}$ (see \S 4.5), we see that 
$\mathbb{G}(b') = \Phi _{\lambda ,w} \bigl( \widetilde{G}_{\Tilde{\lambda}}
(\omega (\Tilde{b})) \bigl) = \Phi _{\lambda ,w}
\bigl( \omega \bigl( \widetilde{G}_{\Tilde{\lambda}}(\Tilde{b}) \bigl) \bigl)
=\Phi _{\lambda ,w}\bigl( \widetilde{G}_{\Tilde{\lambda}}(\Tilde{b}) \bigl)
=\mathbb{G}(\Tilde{b})$. Therefore, we deduce that  
$\mathbb{G}(b') \equiv b \mod qL(\lambda)$.
This completes the proof of Proposition \ref{4.6.1}. \qed

\vspace{2mm}
In the same way as above, we can also show that the map $\Psi _{\lambda , w}$ 
has properties similar to those for $\Phi _{\lambda , w}$, and hence induces the map 
$\overline{\Psi} _{\lambda , w}: B_w (\lambda ) \rightarrow \widetilde{B}_{\Tilde{w}}(\Tilde{\lambda})$.
In fact, this is identical to the map $\widetilde{\ \ }$.
\begin{cor}\label{4.6.2}
With the notation above, the following hold.
\begin{itemize}
\item[(1)]
The maps $\Phi _{\lambda , w}$ and $\Psi _{\lambda , w}$ induce maps 
$\overline{\Phi }_{\lambda , w}: \widetilde{B}_{\Tilde{w}}(\Tilde{\lambda})
\twoheadrightarrow B_w (\lambda )$ and 
$\overline{\Psi } _{\lambda , w}: B_w (\lambda) \hookrightarrow 
\widetilde{B}_{\Tilde{w}}(\Tilde{\lambda})$
such that $\overline{\Phi } _{\lambda , w} \circ \overline{\Psi }_{\lambda , w} 
= \mathrm{id} _{B_w (\lambda )}$.
\item[(2)]
The equalities $G_{\lambda} \circ \overline{\Phi }_{\lambda , w}=
\Phi _{\lambda , w} \circ \widetilde{G}_{\Tilde{\lambda}}$ and 
$\widetilde{G}_{\Tilde{\lambda}} \circ \overline{\Psi } _{\lambda , w}
=\Psi _{\lambda , w} \circ G_{\lambda}$ hold.
\item[(3)]
Let $b \in B_w (\lambda)$, $\Tilde{b}=\overline{\Psi } _{\lambda , w}(b)$, and $\omega \in \Omega _{\lambda ,w}$. 
Then, the equality $\mathbb{G}\bigl( \omega (\Tilde{b}) \bigl) = G_{\lambda}(b)$ holds.
\end{itemize}
\end{cor}

\section{Proofs of main results}
\subsection{Proof of Theorem \ref{thm3}}
\textbf{Proof of Theorem \ref{thm3}.}\ 
If we write the element 
$\Psi _{\lambda ,w} (v) \in \widetilde{V}_{\Tilde{w}}(\Tilde{\lambda })$
for $v\in V_w (\lambda )$ as
$\Psi _{\lambda ,w} (v) = \sum c_{\omega (\Tilde{b})} \widetilde{G}_{\Tilde{\lambda }}
\bigl( \omega (\Tilde{b}) \bigl),\ \text{with}\ c_{\omega (\Tilde{b})} \in \mathbb{C}(q)$, 
then we have 
$v= \Phi _{\lambda ,w}(\Psi _{\lambda ,w} (v))=
\sum c_{\omega (\Tilde{b})}  G_{\lambda }(b)$
by Corollary \ref{4.6.2} (3). Therefore, we conclude that 
$V_{w}(\lambda )=\bigoplus _{b\in B_{w}(\lambda )}\mathbb{C}(q)G_{\lambda }(b).$
This proves Theorem \ref{thm3}. \qed

\subsection{Proof of Theorem \ref{thm4}}
Before starting the proof of Theorem \ref{thm4}, we fix some notation. 
Let $\mathbb{B}_w (\lambda)$, 
$\widetilde{\mathbb{B}}_{\Tilde{w}}(\Tilde{\lambda})$, and 
$\widetilde{\mathbb{B}}_{\Tilde{w}} ^0 (\Tilde{\lambda})$ denote 
the subsets of path crystals corresponding to $B_w (\lambda)$, 
$\widetilde{B}_{\Tilde{w}}(\Tilde{\lambda})$, and 
$\widetilde{B}_{\Tilde{w}} ^0 (\Tilde{\lambda})$ via the 
isomorphisms 
$\mathbb{B}(\lambda ) \cong B(\lambda )$ and 
$\widetilde{\mathbb{B}}(\Tilde{\lambda }) \cong 
\widetilde{B}(\Tilde{\lambda})$ of crystals, respectively.
Define subsets $\mathbb{B}_m$ of $\mathbb{B}_w (\lambda )$ for $m=1,2$, by 
$$\mathbb{B}_1 :=
\bigl \{ F_{w_{k-1}} ^{\mathbf{m}_{k-1}}
f_{i_{k-1}} ^{\epsilon _{k-1}} \cdots F_{w_1} ^{\mathbf{m}_1} 
f_{i_1} ^{\epsilon _1} F_{w_0} ^{\mathbf{m}_0}
\pi _{\lambda } \in \mathbb{B}_w (\lambda ) \ \bigl| \ \mathbf{m}_s \in 
\mathbb{Z}_{\ge 0} ^{\ell (w_s)},\ 0\le \epsilon _s \le a_s 
\ \text{for\ each}\ s \bigl \} ,$$
$$\mathbb{B}_2 := \{ f_{i_k} ^{\epsilon } \pi \in \mathbb{B}_w (\lambda ) \mid 
\pi \in \mathbb{B}_1 ,\ 0\le \epsilon \le a_k \} ,$$
and set $\widetilde{\mathbb{B}} _m := \overline{\Psi}_{\lambda , w}(\mathbb{B}_m)
\subset \widetilde{\mathbb{B}}_{\Tilde{w}} ^0 (\Tilde{\lambda })$, $m=1,2$. \\ \ \\
\textbf{Proof of Theorem \ref{thm4}.}\ 
If $w \in \mathcal{W}_{re}$, then 
$\mathrm{ch}\ \! \mathbb{B}_w (\lambda )= \mathcal{D}_w (e^{\lambda })$ 
by [\textbf{Kas2}, Proposition 3.3.5], and the assertion of Theorem \ref{thm4} follows in this case. 
Let $w = w_k r_{i_k} ^{a_k} \cdots w_1 r_{i_1} ^{a_1} w_0$ be a 
minimal dominant reduced expression. We assume that the equality 
$\mathrm{ch}\ \! \mathbb{B}_1 = 
\mathcal{D}_{w_{k-1}} \mathcal{D}_{i_{k-1}}^{(a_{k-1})}
\cdots \mathcal{D}_{w_1} \mathcal{D}_{i_1 }^{(a_1)}
\mathcal{D}_{w_0} (e^{\lambda })$ holds, and show that 
$\mathrm{ch}\ \! \mathbb{B}_w (\lambda )
=\mathcal{D}_{w_k} \mathcal{D}_{i_k} ^{(a_k)}( \mathrm{ch}\ \! \mathbb{B}_1)$.

First, we show that 
$\mathrm{ch}\ \! \mathbb{B}_2 =\mathcal{D}_{i_k} ^{(a_k)}(\mathrm{ch}\ \! \mathbb{B}_1 )$.
Note that $f_{i_k} ^{\epsilon } \pi \in \mathbb{B}_2$
is not $\mathbf{0}$ for $\epsilon \ge 1$ if and only if 
$\alpha _{i_k}^{\vee}\bigl( \pi (1) \bigl) >0$
(see [\textbf{JL}, Lemma 4.1.6 (1)]). Therefore, we deduce that 
\begin{align}\label{eq11}
\mathrm{ch}\ \! \mathbb{B}_2 = 
\sum _{\begin{subarray}\ \hspace{5mm} \pi \in \mathbb{B}_1 \\
\alpha _{i_k}^{\vee}(\pi (1)) >0 \end{subarray}}
\sum _{\epsilon =0} ^{a_k} e^{f_{i_k} ^{\epsilon } \pi (1)}
+\sum _{\begin{subarray}\ \hspace{5mm} \pi \in \mathbb{B}_1 \\
\alpha _{i_k}^{\vee}(\pi (1)) =0 \end{subarray}} e^{\pi (1)},
\end{align}
and that the right-hand side of (\ref{eq11}) is identical to 
$\mathcal{D}_{i_k} ^{(a_k)} (\mathrm{ch}\ \! \mathbb{B}_1 )$
by the definition of $\mathcal{D}_{i_k} ^{(a_k)}$.

Now, we show that 
$\mathrm{ch}\ \! \mathbb{B}_w (\lambda )
=\mathcal{D}_{w_k} (\mathrm{ch}\ \! \mathbb{B}_2 )$.
By [\textbf{Kas2}, Proposition 3.3.5], we have 
$\mathrm{ch}\ \! \widetilde{\mathbb{B}}_{\Tilde{w}}(\Tilde{\lambda })
=\mathcal{D}_{\Tilde{w}}(e^{\Tilde{\lambda }})
=\mathcal{D}_{\Tilde{w}_k} \mathcal{D}_{\Tilde{w}_k ^{-1} \Tilde{w}}(e^{\Tilde{\lambda }})$.
Then, we can deduce that 
$\mathrm{ch}\ \! \widetilde{\mathbb{B}}_{\Tilde{w}} ^0 (\Tilde{\lambda })
=\mathcal{D}_{\Tilde{w}_k}(\mathrm{ch}\ \! \widetilde{\mathbb{B}}_2 )$
since $\widetilde{\mathbb{B}}_{\Tilde{w}} ^0 (\Tilde{\lambda })
=\{ F_{\Tilde{w}_k } ^{\mathbf{m}} \eta \mid \eta \in 
\widetilde{\mathbb{B}}_2 , \ \mathbf{m} \in \mathbb{Z}_{\ge 0} ^{\ell (\Tilde{w}_k)} \} 
\subset \widetilde{\mathbb{B}}_{\Tilde{w}}(\Tilde{\lambda })$, 
and $\widetilde{\mathbb{B}}_2$ has the \textit{string property}
(see [\textbf{Jo}, Lemma in \S 4.4]).
Therefore, we deduce that
$\mathrm{ch}\ \! \mathbb{B}_w (\lambda )
=\mathrm{ch}\ \! \overline{\Phi}_{\lambda , w}\bigl( \widetilde{\mathbb{B}}_{\Tilde{w}} ^0 (\Tilde{\lambda })\bigl)
=\mathcal{D}_{w_k} (\mathrm{ch}\ \! \overline{\Phi}_{\lambda , w} (\widetilde{\mathbb{B}}_2 ))
=\mathcal{D}_{w_k} (\mathrm{ch}\ \! \mathbb{B}_2 )$.

Consequently, we obtain
$\mathrm{ch}\ \! \mathbb{B}_w (\lambda )
=\mathcal{D}_{w_k} \mathcal{D}_{i_k} ^{(a_k)}(\mathrm{ch}\ \! \mathbb{B}_1 )$.
This proves Theorem \ref{thm4} by induction on $k$. \qed

\end{document}